\documentclass[12pt]{article}
\usepackage{amsfonts}

\usepackage{latexsym}
\usepackage{amssymb}
\usepackage{amsmath}
\usepackage{amsthm}
\usepackage[mathscr]{eucal}
\usepackage{graphicx}
\usepackage{hyperref}
\usepackage{caption}
\usepackage{subcaption}

\newtheorem{Lemma}{Lemma}

\newtheorem{Theorem}[Lemma]{Theorem}

\renewcommand{\qed}{\hfill{\ \ \rule{2mm}{2mm}} \vspace{0.2in}}

\newcommand{\ind}{1\hspace{-2.3mm}{1}}

\begin{document}

\title{Minimum spanning trees across dense cities}
\author{ \textbf{Ghurumuruhan Ganesan}
\thanks{E-Mail: \texttt{gganesan82@gmail.com} } \\
%EndAName
\ \\
New York University, Abu Dhabi }
\date{}
\maketitle

\begin{abstract}
Consider~\(n\) nodes distributed independently across~\(N\) cities
contained with the unit square~\(S\) according to a distribution~\(f.\)
Each city is modelled as an~\(r_n \times r_n\) square contained within~\(S\)
and~\(MSTC_n\) denotes the length of the minimum spanning tree containing
all the~\(n\) nodes. We use approximation methods to 
obtain variance estimates for~\(MSTC_n\)
and prove that if the cities are well-connected and densely populated
in a certain sense, then~\(MSTC_n\) appropriately centred and scaled
converges to zero in probability.

Using the proof techniques,
we alternately derive corresponding results for the length~\(MST_n\)
of the minimum spanning tree
for the usual case when the nodes are independently distributed throughout
the unit square~\(S.\) In particular, we obtain that the variance
of~\(MST_n\) grows at most as a power of the logarithm of~\(n\) and use
a subsequence argument to get almost sure convergence of~\(MST_n\)
appropriately centred and scaled. %\(\frac{1}{\sqrt{n}}(MST_n - \mathbb{E}MST_n)\)

 %and analogous results hold for minimum spanning trees and paths.

%In this paper, we study the structure of left-right crossings
%of the random geometric graph \(G = G(n,r_n)\) of \(n\) nodes
%uniformly distributed in \(S = [0,1]^2\) with \(r_n = \epsilon\sqrt{\frac{\log{n}}{n}}\)
%for some \(\epsilon > 0.\) Tiling \(S\) horizontally and
%vertically into rectangles of length \(1\) and width \(Mr_n,\) we
%show that each rectangle has a left-right crossing of edges with
%high probability if \(M\) is sufficiently large.
%We call the resulting subgraph to be a ``backbone" of \(G.\)

%The techniques we use to construct the backbone has quite a few applications.
%As a first, we show that the diameter of second largest component in \(G\)
%is \(O(1)\) with high probability. Secondly,
\vspace{0.1in} \noindent \textbf{Key words:} Minimum spanning tree, dense cities.

\vspace{0.1in} \noindent \textbf{AMS 2000 Subject Classification:} Primary:
60J10, 60K35; Secondary: 60C05, 62E10, 90B15, 91D30.
\end{abstract}

\bigskip

\setcounter{equation}{0}
\renewcommand\theequation{\thesection.\arabic{equation}}
\section{Introduction}\label{intro}
The study of minimum weight spanning trees of a graph arise in many applications
and many analytical results have been derived regarding the weight of the minimum spanning tree (MST)
for various types of weighted graphs. In this paper, we concern with Euclidean random graphs
where nodes are distributed randomly across the unit square and the goal
is to determine the overall length of the MST. Beardwood et al used subadditive ergodic type
results to obtain that the minimum length of the MST~\(\frac{MST_n}{\sqrt{n}}\) appropriately scaled
converges to a constant a.s.\ as~\(n \rightarrow \infty.\) For more results on MST,
we refer to Steele (1988, 1993), Alexander (1996), Kesten and Lee (1996).

Because of its practical importance, many algorithms have been proposed over the years
to compute the MST for various kinds of graphs. For example, Kruskal's algorithm (Cormen et al (2001))
iteratively adds edges to a sequence of increasing subtree of the original graph
until a spanning tree is obtained. Much of the analytical literature is devoted
to nodes distributed on regular shapes like circles or squares where subadditive
techniques are applicable.

In the first part of this paper, we consider a slightly different problem where nodes are distributed
across small cities distributed throughout the unit square~\(S.\) The cities
are not necessarily regularly placed and therefore subadditive techniques are not directly applicable.
We use approximation methods to obtain sharp bounds for the length
of the minimum spanning tree and thereby deduce the corresponding convergence properties.

%Using our proof techniques, we also obtain variance estimates

\subsection*{Model Description}
\subsubsection*{Structure of the cities}
\begin{figure}[tbp]
\centering
%\fbox{
\includegraphics[width=2in, trim= 100 180 50 100, clip=true]{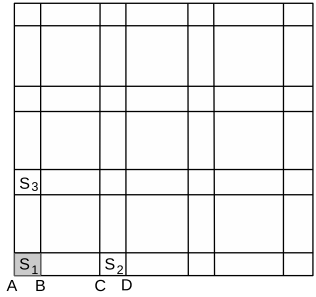}
%}
\caption{Tiling~\(S\) into~\(r_n \times r_n\) squares with an inter-square distance of~\(s_n.\)}
\label{sq_plc}
\end{figure}

For integer~\(n \geq 1,\) let~\(r_n\) and~\(s_n\) be real numbers such that~\(\frac{1-r_n}{r_n+s_n}\)
is an integer. Tile the unit square~\(S\) regularly into~\(r_n \times r_n\) size squares
in such a way that the distance between any two squares is at least~\(s_n\)
as shown in Figure~\ref{sq_plc}. In Figure~\ref{sq_plc}, the grey square is of
size~\(r_n  \times r_n,\)  the segment~\(AB\) has length~\(r_n\) and the segment~\(BC\)
has length~\(s_n.\) The~\(r_n \times r_n\) squares are
called \emph{cities} and the term~\(s_n\) denotes the
\emph{intercity distance}.

Label the~\(r_n \times r_n\) squares (cities) as~\(\{S_l\}\)
and identifying the centres of the squares~\(\{S_l\}\) with vertices in~\(\mathbb{Z}^2,\)
we obtain a corresponding subset of vertices~\(\{z_l\} \subset \mathbb{Z}^2.\)
For example, in Figure~\ref{sq_plc}, identify the centre of the square labelled~\(S_1\)
with~\((0,0),\) the centre of~\(S_2\) with~\((1,0),\) the centre of~\(S_3\) with~\((0,1)\)
and so on. Two vertices~\(z_1 = (x_1,y_1)\) and~\(z_2 = (x_2,y_2)\) are \emph{adjacent}
and connected by an edge if~\(|x_1-x_2| + |y_1 - y_2| = 1.\) %In Figure~\ref{sq_plc}
%for example, the vertices~\(z_1 = (0,0)\) and~\(z_2 = (1,0)\) corresponding respectively to
%the centres of the~\(r_n \times r_n\) squares~\(S_1\) and~\(S_2,\) are adjacent. %and so are~\(S_1\) and~\(S_3.\)

Fix~\(N = N(n)\) cities~\(\{S_{j_1},\ldots,S_{j_N}\}\) and let~\(\{z_{j_1},\ldots,z_{j_N}\}\)
be the vertices in~\(\mathbb{Z}^2\) corresponding to
the centres of~\(\{S_{j_i}\}.\) We say that the cities~\(\{S_{j_1},\ldots,S_{j_N}\}\)
are \emph{well-connected} if the corresponding set of vertices~\(\{z_{j_i}\}\)
form a connected subgraph of~\(\mathbb{Z}^2.\) Henceforth, we
assume that~\(\{S_{j_1},\ldots,S_{j_N}\}\) are well-connected
and without loss of generality denote~\(S_{j_i}\) by~\(S_i\)
for~\(1 \leq i \leq N.\)

\subsubsection*{Nodes in the cities}
Let~\(f\) be any density on the unit square~\(S\) satisfying the following conditions:\\
There are constants~\(0 < \epsilon_1 \leq \epsilon_2 < \infty\) such that
\begin{equation}\label{f_eq}
\epsilon_1 \leq \inf_{x \in S} f(x) \leq \sup_{x \in S} f(x) \leq \epsilon_2
\end{equation}
and
\begin{equation}\label{f_eq_tot}
\int_{x \in S} f(x) dx = 1.
\end{equation}
Define the density~\(g_N(.)\) on the~\(N\) cities~\(\bigcup_{1 \leq i \leq N} S_i\)
as
\begin{equation}\label{gn_def}
g_N(x) = \frac{f(x)}{\int_{\cup_{1 \leq j \leq N}S_j} f_j(x) dx}
\end{equation}
for all~\(x \in \bigcup_{1 \leq j \leq N} S_j.\)

%where~\(vol(A)\) refers to the area of~\(A.\)

%For \(1 \leq i \leq n,\) we define \(X_i\) on the probability space \((S, {\cal B}(S), \mathbb{P}_i),\) where \(\mathbb{P}_i\) denotes the Lebesgue measure on \(S.\)
Let~\(X_1,X_2,\ldots,X_n\) be~\(n\) nodes independently
and identically distributed (i.i.d.)
in the~\(N\) cities~\(\{S_j\}_{1 \leq j \leq N},\) each according to
the density \(g_N.\) Define the vector~\((X_1,\ldots,X_n)\) on the probability
space \((\Omega_X, {\cal F}_X, \mathbb{P}).\) Let~\(K_n = K(X_1,\ldots,X_n)\) be the complete graph whose edges are
obtained by connecting each pair of nodes~\(X_i\) and~\(X_j\) by the straight line segment~\((X_i,X_j)\)
with~\(X_i\) and~\(X_j\) as endvertices. The line segment~\(e_{ij} = (X_i,X_j)\)
is the edge between the nodes~\(X_i\) and\(X_j\)
and~\(d(e_{ij})\) denotes the (Euclidean) length of the edge~\((X_i,X_j).\)

%For an edge~\(e \in K_n\) let~\(l(e)\) denote the length of~\(e.\)

Let~\(Y_1,\ldots,Y_t \subset \{X_k\}\) be~\(t\) distinct nodes.
A path~\({\cal P} = (Y_1,\ldots,Y_t)\) is a subgraph of~\(K_n\)
with vertex set~\(\{Y_{j}\}_{1 \leq j \leq t}\)
and edge set~\(\{(Y_{j},Y_{{j+1}})\}_{1 \leq j \leq t-1}.\)
The nodes~\(Y_1\) and~\(Y_t\) are said to be \emph{connected} by
edges of the path~\({\cal P}.\)
The subgraph~\({\cal C}  = (Y_1,Y_2,\ldots,Y_t,Y_1)\) with
vertex set~\(\{Y_{j}\}_{1 \leq j \leq t}\)
and edge set\\~\(\{(Y_{j},Y_{{j+1}})\}_{1 \leq j \leq t-1} \cup \{(Y_t,Y_1)\}\)
is said to be a \emph{cycle}.

A subgraph~\({\cal T}\) of~\(K_n\)
with vertex set~\(\{Y_i\}_{1 \leq i \leq t}\) and edge set~\(E_{\cal T}\)
is said to be a \emph{tree} if
the following two conditions hold:\\
\((1)\) The graph~\({\cal T}\) is connected; i.e., any two nodes
in~\({\cal T}\) are connected by a path containing only edges in~\(E_{\cal T}.\)\\
\((2)\) The graph~\({\cal T}\) is acyclic; i.e., no subgraph of~\({\cal T}\) is a cycle.\\
The length of the tree~\({\cal T}\) is the sum of the lengths of the edges in~\({\cal T};\) i.e.,
\begin{equation}\label{len_cyc_def}
L({\cal T}) = \sum_{e \in {\cal T}} d(e) = \frac{1}{2} \sum_{i=1}^{t} l(Y_i,{\cal T}),
\end{equation}
where~\(l(Y_i,{\cal T})\) is the sum of lengths of edges in~\({\cal T}\) containing~\(Y_i\)
as an endvertex.

The tree~\({\cal T}\) is said to be a \emph{spanning tree} if~\({\cal T}\)
contains all the~\(n\) nodes\\\(\{X_k\}_{1 \leq k \leq n}.\) Let~\({\cal T}_n\)
be a spanning tree satisfying
\begin{equation}\label{min_weight_tree}
MSTC_n = L({\cal T}_n) := \min_{{\cal T}} L({\cal T}),
\end{equation}
where the minimum is taken over all spanning trees~\({\cal T}.\)
If there is more than one choice for~\({\cal T}_n,\)
choose one according to a deterministic rule. The tree~\({\cal T}_n\)
is defined to the \emph{minimum spanning tree} (MST) with corresponding length~\(MSTC_n.\)

Letting
\begin{equation}\label{bn_def}
b_n := r_n \sqrt{nN},
\end{equation}
we have the following result.
\begin{Theorem}\label{mst_thm}
Suppose~\(r_n,s_n\) and~\(N = N(n)\) satisfy
\begin{equation}\label{N_est}
r_n^2 \geq \frac{M\log{n}}{n}, \frac{n}{N^2} \longrightarrow 0 \text{ and }\frac{N s_n}{b_n} \longrightarrow 0
\end{equation}
as~\(n \rightarrow \infty,\) for some constant~\(M >0.\) If~\(M =M(\epsilon_1,\epsilon_2) > 0\) is large, then
\begin{equation} \label{conv_mst_prob}
\frac{1}{b_n}\left(MSTC_n - \mathbb{E}MSTC_n\right) \longrightarrow 0 \text{ in probability}
\end{equation}
as~\(n \rightarrow \infty.\) In addition, there are positive constants~\(\{\theta_i\}_{1 \leq i \leq 6}\) such that
\begin{equation}\label{exp_mstc}
\theta_1 b_n \leq \mathbb{E}MSTC_n \leq \theta_2 b_n,
\end{equation}
\begin{equation}\label{eq_mst2}
\mathbb{P}\left(MSTC_n \geq  \theta_3 b_n \right) \geq 1 - e^{-\theta_4 N}
\end{equation}
and
\begin{equation}\label{eq_mst3}
\mathbb{P}\left(MSTC_n \leq  \theta_5 b_n \right) \geq 1 - \exp\left(-\theta_6 \frac{n}{N}\right)
\end{equation}
for all~\(n\) large. %and for some constant~\(\theta > 0.\)
\end{Theorem}
In words, if the cities are wide and dense enough, then the centred and scaled minimum length
of the MST converges to zero in probability. %As an example,Condition~(\ref{N_est})
%is satisfied if~\(r_n^2 = \frac{M\log{n}}{n}\)

\subsection*{Unconstrained MST}
There are~\(n\) nodes~\(\{X_k\}_{1 \leq k \leq n}\) independently distributed in the unit square~\(S,\)
each according to the distribution~\(f\) satisfying~(\ref{f_eq}). Let~\({\cal T}_n\) and~\(MST_n\)
denote the minimum spanning tree and its length, respectively, as defined in~(\ref{min_weight_tree}).
Beardwood et al (1959) use subadditive techniques to study the convergence of the
ratio~\(\frac{MST_n}{\sqrt{n}} \longrightarrow \beta\) for some constant~\(\beta > 0,\) a.s.
as~\(n \rightarrow \infty.\) Another approach involves the study of concentration of~\(MST_n\) around its mean
via concentration inequalities~(see Steele (1993)). Here we use the techniques
used in the proof of Theorem~\ref{mst_thm} to obtain the following
result.
\begin{Theorem}\label{var_mst_thm} The variance
\begin{equation}\label{var_mst_est_main}
\mathbb{E}\left(MST_n - \mathbb{E} MST_n\right)^2 \leq C(\log{n})^3
\end{equation}
for some constant~\(C > 0\) and for all~\(n \geq 1\) and
\begin{equation}\label{as_conv_mst}
\frac{1}{\sqrt{n}}\left(MST_n - \mathbb{E} MST_n\right)
\longrightarrow 0 \text{ a.s. }
\end{equation}
as~\(n \rightarrow \infty.\)
There are positive constants~\(\{\theta_i\}_{1 \leq i \leq 3}\) such that
\begin{equation}\label{exp_mst_u}
\theta_1 \sqrt{n} \leq \mathbb{E}MST_n \leq 3 \sqrt{n},
\end{equation}
\begin{equation}\label{eq_mst1_u}
\mathbb{P}\left(MST_n \leq 3\sqrt{n} \right) =1
\end{equation}
and
\begin{equation}\label{eq_mst2_u}
\mathbb{P}\left(MST_n \geq \theta_2 \sqrt{n} \right) \geq 1 - \exp\left(-\frac{\theta_3 n}{\log{n}}\right)
\end{equation}
for all~\(n\) large.

Moreover, if the nodes are uniformly distributed in~\(S,\)
\begin{equation}\label{beta_conv}
\frac{\mathbb{E}MST_n}{\sqrt{n}} \longrightarrow \beta
\end{equation}
as~\(n \rightarrow \infty\) for some constant~\(\beta > 0.\)
\end{Theorem}

The paper is organized as follows. In Section~\ref{prelim}, we state the preliminary estimates
needed for the proofs of the main Theorems. In Section~\ref{pf_mst}, we prove Theorem~\ref{mst_thm}
and in Section~\ref{pf_mst_un}, we prove Theorem~\ref{var_mst_thm}.

%WRT SMALL PROF OF COR +etC...

\setcounter{equation}{0}
\renewcommand\theequation{\thesection.\arabic{equation}}
\section{Preliminary estimates}\label{prelim}
We first derive a deterministic estimate based on the strips method used throughout.

\subsection*{Strips estimate}
Suppose there are~\(a \geq 3\) nodes~\(\{x_i\}_{1 \leq i\leq a}\) placed in a square~\(R\)
of side length~\(b\) such that no two of the nodes share the same~\(x-\)
or~\(y-\)coordinate.  This is a mild condition since if~\(\{X_j\}_{1 \leq j \leq n}\) are i.i.d.\
with density~\(g_N\) as in~(\ref{gn_def}), this condition is satisfied with probability one.
For~\(3 \leq j \leq a\) let~\(K(x_1,\ldots,x_j)\) be the
complete graph with vertex set~\(\{x_i\}_{1 \leq i \leq j}\)
and let~\({\cal T}_{j}\) be a spanning tree of~\(K(x_1,\ldots,x_j)\) such that
\begin{equation}\label{min_tree}
L({\cal T}_j) = \min_{\cal T} L({\cal T}) =: MST(x_1,\ldots,x_j;R),
\end{equation}
where the minimum is taken over all spanning trees of~\(K(x_1,\ldots,x_j)\) and~\(L({\cal T})\)
is the length of the tree~\({\cal T}\) (see~(\ref{len_cyc_def})).

\begin{figure}[tbp]
\centering
%\fbox{
\includegraphics[width=3in, trim= 20 180 50 110, clip=true]{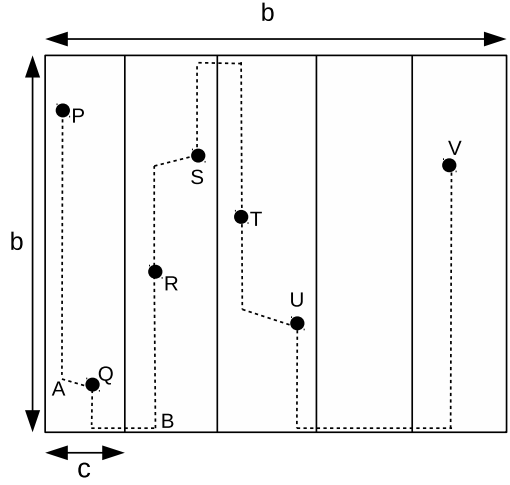}
%}
\caption{Estimating minimum length using strips counting.}
\label{stp_fig}
\end{figure}

We have that
\begin{equation}\label{mst_ab}
MST(x_1,\ldots,x_a;R) \leq 3b \sqrt{a}.
\end{equation}
\emph{Proof of~(\ref{mst_ab})}: Divide the square~\(R\) into vertical rectangles (strips) each of
size~\(c \times b\) so that the number of strips
is~\(\frac{b}{c}\) as shown in Figure~\ref{stp_fig}. Here~\(a= 7\) and without loss of generality
suppose that~\(P,Q, R, S,T,U\) and~\(V,\) are the nodes~\(x_1,x_2,x_3,x_4,x_5,x_6\) and~\(x_7,\) respectively.
The dotted line corresponds to a path
containing all the nodes~\(P,Q,R,S\) and~\(T.\) Starting from the top most node~\(P\) in the first strip,
vertically down in the strip and each time we are close to a node, we ``reach" for the node by a slightly
inclined line.
In Figure~\ref{stp_fig}, the vertical dotted line~\(PA\) is joined to the node~\(Q\) by the inclined
line~\(AQ.\)

Continue vertically down from~\(Q\) until we reach close to the bottom of the strip.
Proceed along a horizontal line until we are directly below the lowest node
in the second strip. In Figure~\ref{stp_fig}, the point~\(B\) is directly below the node~\(R.\)
Continue vertically from~\(B,\) pass through~\(R\) until we reach close to the next node~\(S.\)
Join to~\(S\) by a slightly inclined line and continue this procedure until all nodes
in all strips have been exhausted.

The number of strips is~\(\frac{b}{c}\) and the
sum of the lengths of the vertical lines of~\({\cal P}\) in a particular strip is at most the height of the strip~\(b.\)
Therefore the total length of vertical lines in~\({\cal P}\) is at most~\(\frac{b}{c}b.\)

The total length of the horizontal lines in~\({\cal P}\) is at most~\(b.\)
Finally, each inclined line in~\({\cal P}\) has length at most~\(\frac{c}{\sqrt{2}},\) since the corresponding slope
is at most~\(45\) degrees. Each of the~\(a\) nodes is attached to at most one inclined line
and so the total length of the inclined lines in~\({\cal P}\)
is at most~\(\frac{a c}{\sqrt{2}}.\)

Summarizing, the total length of edges in~\({\cal P}\) is at most~\(\frac{b^2}{c}  + \frac{ac}{\sqrt{2}} + b.\)
By construction, the path~\({\cal P}\)
encounters the nodes~\(x_1,\ldots,x_a\) in that order and so applying triangle inequality
as before, the path~\({\cal P}_0 = (x_1,x_2,\ldots,x_a)\)
with edges being the straight lines~\((x_1,x_{2}), (x_2,x_3),\ldots,(x_a,x_1),\)
has total length no more than the sum of length of edges in~\({\cal P}.\) Thus
\begin{equation}\label{temp_est_strip}
MST(x_1,\ldots,x_a; R) \leq L({\cal P}_0) \leq \frac{b^2}{c}  + \frac{ac}{\sqrt{2}} + b.
\end{equation}
Setting~\(c = \frac{b}{\sqrt{a}}\) in~(\ref{temp_est_strip}), we get that~\(MST(x_1,\ldots,x_a;R)\)
is bounded above by~\(b\sqrt{a} + \frac{b\sqrt{a}}{\sqrt{2}} + b \leq 3b\sqrt{a}, \)
since~\(a \geq 1.\)~\(\qed\)

%We also need another estimate that merges trees within small squares in~\(S.\)

\subsection*{Length of MST within cities}
Recall from discussion prior to~(\ref{N_est}) that~\(n \geq 1\) nodes~\(\{X_k\}_{1 \leq k \leq n}\)
are distributed across the~\(r_n \times r_n\) squares~\(\{S_j\}_{1 \leq j \leq N}\)
according to a Binomial process with intensity~\(g_N\) as defined in~(\ref{gn_def}).
In this subsection, we obtain estimates for the length~\(R_l\)
of the MST containing all the nodes of the square~\(S_l.\)

If~\(p_l\) denotes the probability that a node of~\(\{X_j\}\) occurs inside~\(S_l,\)
then
\begin{equation}\label{pl_def}
\frac{\eta_1}{N} \leq p_l := \frac{\int_{S_l} f(x) dx}{\int_{\cup_j S_j} f(x)  dx} \leq \frac{\eta_2}{N},
\end{equation}
where~\(\eta_1 = \frac{\epsilon_1}{\epsilon_2} \leq \frac{\epsilon_2}{\epsilon_1} = \eta_2\) (see~(\ref{f_eq})).
Therefore if
\begin{equation}\label{nl_def}
N_l = \sum_{i=1}^{n} \ind(X_i \in S_l)
\end{equation}
denotes the number of nodes of~\(\{X_j\}\) in the square~\(S_l,\) then~\(N_l\)
is Binomially distributed with parameters~\(n\) and~\(p_l;\) i.e., for any~\(1 \leq  k \leq n,\)
\begin{equation}\label{bin_dist}
\mathbb{P}(N_l = k) = B(k;n,p_{l}) := {n \choose k} p_l^{k} (1-p_l)^{n-k},
\end{equation}
where~\({n \choose k} = \frac{n!}{k!(n-k)!}\) is the Binomial coefficient.
Moreover,
\begin{equation}\label{nl_exp}
\frac{\eta_1 n}{N} \leq \mathbb{E} N_l = np_l \leq \frac{\eta_2 n}{N}
\end{equation}
by~(\ref{pl_def}).

Let~\(\{Y_j\}_{1 \leq j \leq N_l}\)
be the nodes of~\(\{X_j\}\) present in the square~\(S_l.\)
Formally, if~\(N_l = 0,\) set~\(\{Y_j\}_{1 \leq j \leq N_l} := \emptyset.\)
If~\(N_l \geq 1,\) define~\(N_l\) indices~\(j_1,\ldots,j_{N_l}\) as follows.
Let~\[j_1 = j_1(X_1,\ldots,X_n) := \min\{1 \leq k \leq n :  X_k \in S_l\}\]
be the least indexed node of~\(\{X_k\}\) present in~\(S_l.\)
Let~\[j_2 = \min\{j_1+1 \leq k \leq n : X_k \in S_l\}\] be the next least indexed node of~\(\{X_k\}\)
present in~\(S_l\) and so on. Set~\(Y_i = X_{j_i}\) for~\(1 \leq j \leq N_l.\)

Set~\(R_l = 0\) if~\(N_l \leq 2\) and if~\(N_l \geq 3\) set
\begin{equation}\label{rl_def}
R_l := MST(Y_1,\ldots,Y_{N_l}; S_l)
\end{equation}
where~\(MST(.;.)\) is as defined in~(\ref{min_tree}). The following is the main lemma proved in this subsection.
\begin{Lemma} \label{rl_lemma}
If~\(M > 0\) is arbitrary and~(\ref{N_est}) holds, the following is true: There are positive constants~\(\{\delta_i\}_{1 \leq i \leq 3}\) such that
for all~\(n \geq 2\) and for any~\(1 \leq l \leq N,\)
\begin{equation}\label{del_tn_b}
\delta_1 r_n\sqrt{\frac{n}{N}}\leq \mathbb{E} R_l \leq \delta_2 r_n \sqrt{\frac{n}{N}} \;\;\; \text{ and }
\;\;\; \mathbb{E} R_l^2 \leq \delta_3 \left(r_n \sqrt{\frac{n}{N}}\right)^2.
\end{equation}
Moreover, if
\begin{equation}\label{ul_def}
U_l = U_l(n) := \left\{\frac{\eta_1 n}{2N} \leq N_l \leq \frac{2\eta_2 n}{N}\right\},
\end{equation}
where~\(\eta_1\) and~\(\eta_2\) are as in~(\ref{pl_def}),
then there are positive constants~\(\{\delta_i\}_{i=4,5}\) such that
for all~\(n \geq 2\) and for any~\(1 \leq l \leq N,\)
\begin{equation}\label{ul_est}
\mathbb{P}(U_l) \geq 1- \exp\left(-\delta_4 \frac{n}{N}\right) \text{ and } R_l\ind(U_l) \leq \delta_5 r_n\sqrt{\frac{n}{N}}.
\end{equation}
\end{Lemma}

To prove~ the above Lemma, we perform some preliminary computations.
We first derive bounds for the total number of squares~\(N.\) From~(\ref{N_est}) we have
that~\(r_n^2 \geq \frac{M\log{n}}{n}\) and since all the~\(r_n \times r_n\) squares~\(\{S_l\}_{1 \leq l \leq N}\)
are contained within the unit square~\(S,\) we also have~\(Nr_n^2 \leq 1\) and therefore~\(N \leq \frac{n}{M\log{n}}.\) Similarly from~(\ref{N_est})
we also have that~\(\frac{n}{N^2} \longrightarrow 0\) as~\(n \rightarrow \infty\) and so~\(N \geq \sqrt{n}\) for all~\(n\) large. Combining we get
\begin{equation}\label{n_N}
\sqrt{n} \leq N \leq \frac{n}{M\log{n}} \text{ and } \frac{n}{N} \geq M\log{n}
\end{equation}
for all~\(n\) large.

For~\(k \geq 2,\) let~\(D_l(k)\) be the expected minimum distance between the node~\(Y_k\)
and every other node in~\(S_l,\) given that there are~\(N_l = k\) nodes in~\(S_l;\) i.e.,
\begin{equation}\label{pf_exp1}
D(k) = D_{l}(k) := \mathbb{E}\left(d(Y_{k},\{Y_u\}_{1 \leq u \leq k-1}) | N_l = k\right),
\end{equation}
where~\(d(A,B) = \min_{x \in A, y \in B} d(x,y)\) is the minimum distance between finite sets~\(A\) and~\(B.\)
We have the following properties.\\
\((b1)\) For any~\(k \geq 2\) and~\(1 \leq l \leq N,\) the term
\begin{equation}\label{d_low}
D_l(k) \geq  \int_0^{\frac{r_n}{\sqrt{\delta}}} \left(1-\pi \eta_2 \left(\frac{r}{r_n}\right)^2\right)^{k-1} dr
\end{equation}
where~\(\eta_2 = \frac{\epsilon_2}{\epsilon_1}\) is as in~(\ref{pl_def}).\\\\
\((b2)\) There are positive constants~\(\gamma_i,1 \leq i \leq 3\) such that for any~\(k \geq 2\) and~\(1 \leq l \leq N,\) the minimum distance
\begin{equation}\label{x2_est}
\gamma_1 \frac{r_n}{\sqrt{k}} \leq D_l(k) \leq \gamma_2 \frac{r_n}{\sqrt{k}}
\text{ and } \mathbb{E}\left(d^2(Y_{k},\{Y_u\}_{1 \leq u \leq k-1}) | N_l = k\right) \leq \gamma_3 \frac{r_n^2}{k}.
\end{equation}
The proof of~\((b1)-(b2)\) uses the fact that given~\(N_l = k,\) the nodes in~\(S_l\)
are independently distributed in~\(S_l\) with distribution~\(f;\) i.e.,
\begin{equation}\label{dl_zeq}
D_l(k) = \mathbb{E}d(Z_k,\{Z_j\}_{1 \leq j \leq k-1})
\end{equation}
where~\(\{Z_i\}_{1 \leq i \leq k}\) are i.i.d.\
with distribution
\begin{equation}\label{z_dist}
\mathbb{P}(Z_1 \in A) = \frac{\int_{A \cap S_l} f(x) dx}{\int_{S_l} f(x) dx}.
\end{equation}

%This is possible since~\(\frac{n}{N} \longrightarrow \infty\) as~\(n \rightarrow \infty\)~(see~(\ref{n_N})). %In~(\ref{ind_tl_b}), the term~\(D_l(k)\) is the expected minimum distance as defined in~(\ref{pf_exp1}).

%\begin{equation}\label{ind_tl2_b}
%\mathbb{E} T_l  \geq C r_n\sqrt{\frac{n}{N}}
%\end{equation}
%for all~\(n\) large and some constant~\(C >0.\)

Use Fubini's theorem and~(\ref{z_dist}) to write
\begin{equation}\label{dl_alt}
D_l(k) = \frac{1}{\int_{S_l} f(x) dx} \int_{S_l} \mathbb{E}d(x,\{Z_j\}_{1 \leq j \leq k-1})f(x)dx,
\end{equation}
where~\(\mathbb{E}d(x,\{Z_j\}_{1 \leq j \leq k-1}) = \int_{0}^{\infty} \mathbb{P}(d(x,\{Z_j\}_{1 \leq j \leq k-1}) \geq r) dr.\)
For any~\(x \in S_l,\) the minimum distance from~\(x\) to~\(\{Z_1,\ldots,Z_{k-1}\}\) is at least~\(r\) if and
only if~\(B(x,r) \cap S_l\) contains no point of~\(\{Z_j\}_{1 \leq j \leq k-1}.\) Here~\(B(x,r)\) is the
ball of radius~\(r\) centred at~\(x.\) Wherever the point~\(x \in S_l,\) the area of~\(B(x,r) \cap S_l\)
is at most~\(\pi r^2\) and so together with~(\ref{f_eq}), we then get that
\begin{equation}
\mathbb{P}(d(x,\{Z_j\}_{1 \leq j \leq k-1}) \geq r) = \left(1-\frac{\int_{B(x,r) \cap S_l} f(x)dx}{\int_{S_l} f(x)dx}\right)^{k-1}
\label{dx_min}
\nonumber
\end{equation}
is bounded below by~\(\left(1-\pi \eta_2\frac{r^2}{r_n^2}\right)^{k-1},\)
where~\(\eta_2 =\frac{\epsilon_2}{\epsilon_1}\) is as in~(\ref{pl_def}). This proves~(\ref{d_low}).

To prove the lower bound for~\(D_l(k)\) in~(\ref{x2_est}) of~\((b2),\) fix \(k \geq 2\)  and use~(\ref{d_low}) to get that
\begin{equation}
D_l(k) \geq \int_{0}^{\frac{r_n}{\sqrt{\delta k}}}\left(1-\delta\left(\frac{r}{r_n}\right)^2\right)^{k-1} dr
\geq \int_{0}^{\frac{r_n}{\sqrt{\delta k}}}\left(1-\frac{1}{k}\right)^{k-1} dr \geq \frac{e^{-1} r_n}{\sqrt{\delta k}} \nonumber`
\end{equation}
for all~\(n\) large. The final estimate is obtained by using~\(\left(1-\frac{1}{r}\right)^{r-1} \geq e^{-1}\) for all~\(r \geq 2.\)

For the upper bound for~\(D_l(k)\) in~(\ref{x2_est}),  again use~(\ref{dx_min}) and the fact that~\(B(x,r) \cap S_l\) has area at least~\(\frac{\pi r^2}{4}\)
no matter where the position of~\(x,\) to get
\[\mathbb{P}(d(x,\{Z_j\}_{1 \leq j \leq k-1}) \geq r) \leq \left(1-\frac{\pi}{4\epsilon_1} \left(\frac{r}{r_n}\right)^2\right)^{k-1} \leq \exp\left(-\frac{\pi(k-1)}{4\epsilon_1 r_n^2} r^2\right)\] and so~\(D_l(k) \leq \int_{0}^{\infty} \exp\left(-\frac{\pi(k-1)}{4\epsilon_1 r_n^2} r^2\right) dr \leq \frac{C r_n}{\sqrt{k-1}} \leq \frac{2C r_n}{\sqrt{k}}\)
for all~\(k \geq 2\) and for some positive constant~\(C,\) not depending on~\(k\) or~\(l.\)

Finally for the second moment estimate in~(\ref{x2_est}), we argue analogous to~(\ref{pf_exp1}) and get that the
term~\(\mathbb{E}\left(d^2(Y_{k},\{Y_u\}_{1 \leq u \leq k-1}) | N_l = k\right) \) equals
\begin{equation}\label{pf_exp2}
\mathbb{E}d^2(Z_k,\{Z_j\}_{1 \leq j \leq k-1}) = \frac{1}{\int_{S_l} f(x) dx} \int_{S_l} \mathbb{E}d^2(x,\{Z_j\}_{1 \leq j \leq k-1}) f(x) dx
\end{equation}
where~\(\{Z_i\}_{1 \leq i \leq k}\) are i.i.d.\ with distribution as in~(\ref{z_dist}).
Arguing as in the previous paragraph we get that
\begin{equation}
\mathbb{E}(d^2(x,\{Z_j\}_{1 \leq j \leq k-1})) = \int r \mathbb{P}(d(x,\{Z_j\}_{1 \leq j \leq k-1}) \geq r ) dr \nonumber
\end{equation}
is bounded above by~\(\int_{0}^{\infty} r\exp\left(-\frac{\pi(k-1)}{4\epsilon_1 r_n^2} r^2\right) dr \leq \frac{C r_n^2}{k}\)
for some positive constant~\(C,\) not depending on~\(k\) or~\(x.\) This proves
the desired bound for the second moment in~(\ref{x2_est}).~\(\qed\)

\emph{Proof of Lemma~\ref{rl_lemma}}: The proof of the first estimate in~(\ref{ul_est}) follows from standard Binomial estimates and the estimate for~\(\mathbb{E}N_l\) in~(\ref{nl_exp}) (see Corollary~A.1.14, pp. 312, Alon and Spencer (2008)). The proof the second estimate in~(\ref{ul_est}) follows from the strips estimate~(\ref{mst_ab}) with~\(a = \frac{2\eta_2 n}{N}\) and~\(b = r_n.\)

%and~\(c = r_n\sqrt{\frac{N}{n}}.\)

To prove the first estimate of~(\ref{del_tn_b}) assume~\(N_l \geq 3\) and
recall that~\(\{Y_u\}_{1 \leq u \leq N_l}\) are the nodes of the Binomial process in the square~\(S_l\) (see paragraph prior to~(\ref{pf_exp1})).
Let~\({\cal R}_l\) denote the MST of length~\(R_l\)
containing the nodes~\(\{Y_u\}_{1 \leq u \leq N_l}.\)  If~\(l(Y_u,{\cal R}_l),1 \leq u \leq N_l\) is the sum of length of the edges containing~\(Y_u\)
as an endvertex then~\(l(Y_u,{\cal R}_l) \geq d(Y_u,\{Y_v\}_{v \neq u}),\) the minimum distance of~\(Y_u\) from all
the other nodes in~\(S_l\) as defined in~(\ref{pf_exp1}).

From~(\ref{len_cyc_def}),~\(R_l = L({\cal R}_l) = \frac{1}{2} \left(\sum_{u=1}^{N_l} l(Y_u,{\cal R}_l) \right) \geq \frac{1}{2}\left(\sum_{u=1}^{N_l}d(Y_u,\{Y_v\}_{v  \neq u })\right)\)
and so
\begin{equation}
\mathbb{E} R_l  = \sum_{k \geq 2} \mathbb{E} R_l \ind(N_l = k) \geq \frac{1}{2}\mathbb{E} \sum_{k \geq 2} \sum_{u=1}^{k} d(Y_u,\{Y_v\}_{v  \neq u })\ind(N_l = k). \label{ind_tl_b2}
\end{equation}
Recalling the definition of~\(D_l(k)\) in~(\ref{pf_exp1}) we then get
\begin{equation}
\mathbb{E} R_l  \geq \frac{1}{2} \sum_{k \geq 2} \mathbb{P}(N_l = k) kD_l(k) \geq \frac{1}{2}\sum_{\frac{\eta_1 n}{2N} \leq k \leq \frac{2\eta_2 n}{N}}\mathbb{P}(N_l = k) kD_l(k), \label{ind_tl_b}
\end{equation}
provided~\(n\) is large enough so that~\(\frac{\eta_1 n}{2N} \geq \frac{\eta_1}{2} M\log{n} \geq 2,\) the middle estimate being true because of~(\ref{n_N}). %This is possible since~\(\frac{n}{N} \longrightarrow \infty\) as~\(n \rightarrow \infty\)~(see~(\ref{n_N})). %In~(\ref{ind_tl_b}), the term~\(D_l(k)\) is the expected minimum distance as defined in~(\ref{pf_exp1}).

Using the estimate~\(D_l(k) \geq \frac{\gamma_1 r_n}{\sqrt{k}}\) (see~(\ref{x2_est})) in~(\ref{ind_tl_b}) we then get
that~\(\mathbb{E}R_l\) is bounded below by
\[\gamma_1 r_n\sum_{\frac{\eta_1 n}{2N} \leq k \leq \frac{2\eta_2 n}{N}}\mathbb{P}(N_l = k) \sqrt{k}  \geq \gamma_1 r_n \sqrt{\frac{\eta_1 n}{2N}} \sum_{\frac{\eta_1 n}{2N} \leq k \leq \frac{2\eta_2 n}{N}}\mathbb{P}(N_l = k),\]
which in turn is bounded below by~\(\gamma_1 r_n \sqrt{\frac{\eta_1 n}{2N}}  \left(1-e^{-C\frac{n}{N}}\right)\) for some constant~\(C > 0,\) by~(\ref{ul_est}). Since~\(\frac{n}{N}  \longrightarrow \infty\) as~\(n \rightarrow \infty,\) (see~(\ref{n_N})), this proves the lower bound for~\(\mathbb{E} R_l\) in~(\ref{del_tn_b}).

%\begin{equation}\label{ind_tl2_b}
%\mathbb{E} T_l  \geq C r_n\sqrt{\frac{n}{N}}
%\end{equation}
%for all~\(n\) large and some constant~\(C >0.\)

To prove the upper bound of~\(\mathbb{E}R_l\) in~(\ref{del_tn_b}), we argue as follows. If the number of nodes~\(N_l \leq \frac{2\eta_2 n}{N},\) then from~(\ref{ul_est}),~\(R_l \leq C r_n \sqrt{\frac{n}{N}}\) for some constant~\(C > 0.\)
If~\(N_l \geq \frac{2\eta_2 n}{N},\) then~\(R_l\leq N_l r_n \sqrt{2},\) since there are at most~\(N_l-1 \leq N_l\)
edges in the MST~\({\cal R}_l\) of length~\(R_l\) and each such edge
has both endvertices in the~\(r_n \times r_n\) square~\(S_l\) and therefore has length at most~\(r_n\sqrt{2}.\)
Thus
\begin{equation}\label{tij_est_temp_b}
\mathbb{E} R_l \leq C r_n \sqrt{\frac{n}{N}} + r_n\sqrt{2} \mathbb{E} \left(N_l\ind\left(N_l > \frac{2\eta_2 n}{N}\right)\right)
\leq C r_n \sqrt{\frac{n}{N}} + r_n\sqrt{2} \mathbb{E} (N_l\ind(U_l^c)),
\end{equation}
where~\(U_l\) is as defined in~(\ref{ul_def}).

Recall from discussion following~(\ref{nl_def}) that~\(N_l\) is Binomially distributed with parameters~\(n\) and~\(p_l\) and so by standard Binomial estimates~\(\mathbb{E}N_l^2 \leq C(np_l)^2 \leq \frac{Cn^2}{N^2}\) for some constant~\(C > 0,\) by~(\ref{pl_def}).
Using Cauchy-Schwarz inequality and the estimate for~\(\mathbb{P}(U_l)\) in~(\ref{ul_est}), we therefore get
\begin{equation} \label{ul_2233}
\mathbb{E} N_l \ind(U_l^c) \leq \left(\mathbb{E}N_l^2\right)^{\frac{1}{2}} \left(\mathbb{P}(U_l^c)\right)^{\frac{1}{2}} \leq C_1\frac{n}{N}\exp\left(-C_2\frac{n}{N}\right) \leq \sqrt{\frac{n}{N}},
\end{equation}
for all~\(n\) large and for some positive constants~\(C_1,C_2.\) The final inequality in~(\ref{ul_2233}) is true
since~\(\frac{n}{N} \longrightarrow \infty\) as~\(n \rightarrow \infty\) (see~(\ref{n_N})).
Substituting~(\ref{ul_2233}) into~(\ref{tij_est_temp_b}) gives the upper bound for~\(\mathbb{E}R_l\) in~(\ref{del_tn_b}). The proof of the bound for~\(\mathbb{E}R_l^2\) is analogous as above.~\(\qed\)

Define the covariance between~\(R_{l_1}\) and~\(R_{l_2}\) for distinct~\(l_1\) and~\(l_2\) as
\begin{equation}\label{cov_def}
cov(R_{l_1},R_{l_2}) = \mathbb{E}R_{l_1} R_{l_2} - \mathbb{E} R_{l_1} \mathbb{E}R_{l_2}.
\end{equation}
We need the following result for future use. Recall the constants~\(\epsilon_1,\epsilon_2\) in~(\ref{f_eq}).
\begin{Lemma}\label{cov_lemma}
There is a positive constant~\(M_0 = M_0(\epsilon_1,\epsilon_2)\) large so that the following holds if~(\ref{N_est})
is satisfied with~\(M > M_0:\) There are positive constants~\(C_1,C_2\) such that for all~\(n \geq 2\) and for any~\(1 \leq l_1 \neq l_2 \leq N,\)
\begin{equation}\label{cov_tl_est}
|cov(R_{l_1},R_{l_2})| \leq C_1 \left(\mathbb{E}R_{l_1} R_{l_2}\right)\frac{n}{N^2} \leq C_2 \frac{r_n^2 n^2}{N^3}.
\end{equation}
\end{Lemma}
To prove Lemma~\ref{cov_lemma}, we use Poissonization described in the next subsection.

\subsection*{Poissonization}
Recall from discussion prior to~(\ref{N_est}) that~\(n \geq 1\) nodes~\(\{X_k\}_{1 \leq k \leq n}\)
are distributed  across the~\(r_n \times r_n\) squares~\(\{S_j\}_{1 \leq j \leq N}\)
according to a Binomial process
with intensity~\(g_N(.)\) as defined in~(\ref{gn_def}).
Throughout, we use Poissonization as a tool
to obtain estimates for probabilities of events for the corresponding Binomial process.
We make precise the notions in this subsection.

Let~\({\cal P}\) be a Poisson process
on the squares~\(\cup_{j=1}^{N} S_j\) with intensity function~\(ng_N(.)\)
defined on the probability space~\((\Omega_0,{\cal F}_0, \mathbb{P}_0).\)
If~\(N^{(P)}_l\) be the number of nodes of~\({\cal P}\) present in the square~\(S_l, 1 \leq l \leq N,\)
then
\begin{equation}\label{poi_dist}
\mathbb{P}_0(N^{(P)}_l = k) = Poi(k;np_l) := e^{-np_l}\frac{(np_l)^{k}}{k!},
\end{equation}
where~\(p_l\) is as defined in~(\ref{pl_def}).
Moreover,
\begin{equation}\label{nl_exp_p}
\frac{\eta_1 n}{N} \leq \mathbb{E}_0 N^{(P)}_l = np_l \leq \frac{\eta_2 n}{N}
\end{equation}
by~(\ref{pl_def}).

Let~\(\{Y_j\}_{1 \leq j \leq N^{(P)}_l}\)
be the nodes of~\({\cal P}\) present in the square~\(S_l.\)
Analogous to~(\ref{rl_def}), set~\(R^{(P)}_l = 0\) if~\(N^{(P)}_l \leq 2\) and if~\(N^{(P)}_l \geq 3\) set
\begin{equation}\label{rl_def_p}
R^{(P)}_l := MST(Y_1,\ldots,Y_{N^{(P)}_l}; S_l)
\end{equation}
where~\(MST(.;.)\) is as defined in~(\ref{min_tree}).
The following result is analogous to Lemma~\ref{rl_lemma}.
\begin{Lemma} \label{rl_lemma_poiss}
If~\(M > 0\) is arbitrary and~(\ref{N_est}) holds, the following is true: There are positive constants~\(\{\delta_i\}_{1 \leq i \leq 5}\) such that
for all~\(n \geq 2\) and for any~\(1 \leq l \leq N,\)
\begin{equation}\label{del_tn}
\delta_1 r_n\sqrt{\frac{n}{N}}\leq \mathbb{E}_0 R^{(P)}_l \leq \delta_2 r_n \sqrt{\frac{n}{N}}, \;\;\;\;\;\;\;\;\mathbb{E}_0 \left(R^{(P)}_l\right)^2 \leq \delta_3 \left(r_n \sqrt{\frac{n}{N}}\right)^2
\end{equation}
and
\begin{equation}\label{tn_prob_est}
\mathbb{P}_0\left(R^{(P)}_l \geq \delta_4 r_n \sqrt{\frac{n}{N}}\right) \geq \delta_5.
\end{equation}
\end{Lemma}

%Define the distance function~\(d_B(u,N_l)\) for~\(u \geq 1\) as follows.
%Set~\(d_B(u,N_l) = 0\) for all~\(u \geq N_l+1.\) If~\(N_l = 1,\) set~\(d_B(u,N_l)  = 0\) and
%if~\(N_l \geq 2,\) define~\[d_B(u,N_l)  = d(Y_u,\{Y_v\}_{1 \leq v \leq N_l, v \neq u})\] to be
%the minimum distance of node~\(Y_u\) from all the other nodes of~\(\{X_i\}\) in~\(S_l.\)

\emph{Proof of Lemma~\ref{rl_lemma_poiss}}: The proof of~(\ref{del_tn}) is analogous as in the Binomial case and proceeds as follows.
Define
\begin{equation}\label{ul_def_p}
U^{(P)}_l = U^{(P)}_l(n) := \left\{\frac{\eta_1 n}{2N} \leq N^{(P)}_l \leq \frac{2\eta_2 n}{N}\right\},
\end{equation}
where~\(\eta_1\) and~\(\eta_2\) are as in~(\ref{pl_def}). Analogous to~(\ref{ul_est}), the following bound
is obtained from standard Poisson distribution estimates (see Theorem~A.1.15, pp. 313, Alon and Spencer (2008)):
There is a positive constant~\(\gamma\) such that
for all~\(n \geq 2\) and for any~\(1 \leq l \leq N,\)
\begin{equation}\label{ul_est_p}
\mathbb{P}_0\left(U^{(P)}_l\right) \geq 1- \exp\left(-\gamma \frac{n}{N}\right).
\end{equation}

As in the Binomial case, given~\(N^{(P)}_l = k,\) the nodes of~\({\cal P}\) are i.i.d.\ distributed
according to distribution~(\ref{z_dist}). Therefore for~\(k \geq 2\)
we let~\[D^{(P)}_l(k) = \mathbb{E}_0\left(d(Y_k,\{Y_j\}_{1 \leq j \leq k-1})|N^{(P)}_l  = k\right)\]
and as in~(\ref{pf_exp1}) obtain that
\begin{equation}\label{dlp_k}
D^{(P)}_l(k) = \mathbb{E}d(Z_k,\{Z_j\}_{1 \leq j \leq k-1}) = D_l(k),
\end{equation}
where~\(D_l(k)\) is as defined in~(\ref{pf_exp1}), the random variables~\(\{Z_j\}_{1 \leq j \leq k}\)
are i.i.d.\ with distribution~(\ref{z_dist}) and the final equality in~(\ref{dlp_k}) is true
because of~(\ref{dl_zeq}). Consequently~\(D^{(P)}_l(k)\) also satisfies properties~\((b1)-(b2)\)
and the rest of the proof of~(\ref{del_tn}) is analogous to the Binomial case.

Finally, the estimate in~(\ref{tn_prob_est}) is obtained by using~(\ref{del_tn}) and the Paley-Zygmund inequality
\begin{equation}\label{paley}
\mathbb{P}_0\left(R^{(P)}_l \geq \lambda \mathbb{E}_0 R^{(P)}_l\right)
\geq (1-\lambda)^2\frac{(\mathbb{E}_0R^{(P)}_l)^2}{\mathbb{E}_0\left(R^{(P)}_l\right)^2}
\end{equation} for~\(0 < \lambda <1.\)~\(\qed\)

%In the end, we convert the estimates to the Binomial process.

We now use Poissonization
and obtain intermediate estimates needed to prove Lemma~\ref{cov_lemma}.
Recall from~(\ref{rl_def}) and~(\ref{rl_def_p}) that~\(R_l\) and~\(R^{(P)}_l\) are the lengths of the MSTs containing
all the nodes in the~\(r_n \times r_n\) square~\(S_l, 1 \leq l \leq N\) in the Binomial and the Poisson process, respectively.
\begin{Lemma}\label{poi_bin_diff_lemm}
There is a positive constant~\(M_0 = M_0(\epsilon_1,\epsilon_2)\) large so that the following holds if~(\ref{N_est})
is satisfied with~\(M > M_0:\) There are positive constants~\(C_0,C_1\) and~\(C_2\)
not depending on~\(l\) such that the following estimates hold for all~\(n \geq C_0:\)
For~\(1 \leq l \leq N,\)
\begin{equation}\label{tl_diff}
|\mathbb{E}R_l - \mathbb{E}_0 R^{(P)}_l| \leq C_1\left(\mathbb{E}R_l\right)\left(\frac{n}{N^2}\right)
\leq C_2\left(\frac{r_n n^{3/2}}{N^{5/2}}\right).
\end{equation}
For any~\(1 \leq l_1 \neq l_2 \leq N\)
\begin{equation}\label{tl_diff2}
|\mathbb{E}(R_{l_1}R_{l_2}) - \mathbb{E}_0 (R^{(P)}_{l_1} R^{(P)}_{l_2})| \leq C_1\left(\mathbb{E}R_{l_1}\mathbb{E}R_{l_2}\right)\left(\frac{n}{N^2}\right)
\leq C_2\left(\frac{r_n^2 n^2}{N^3}\right).
\end{equation}
\end{Lemma}

To prove Lemma~\ref{poi_bin_diff_lemm}, we need estimates on the difference between Binomial and Poisson distributions.
For~\(k,l \geq 1\) recall the Binomial distribution~\(B(k; n,p_l)\)
and the Poisson distribution~\(Poi(k; np_l)\)
as defined in~(\ref{bin_dist}) and~(\ref{poi_dist}), respectively.
For~\(k_1,k_2,l_1,l_2 \geq 1,\) let
\begin{equation}\label{bin_dist2}
B(k_1,k_2;n,p_{l_1},p_{l_2}) := {n \choose k_1,k_2} p_{l_1}^{k_1} p_{l_2}^{k_2}(1-p_{l_1} - p_{l_2})^{n-k_1-k_2},
\end{equation}
where~\({n \choose k_1,k_2} = \frac{n!}{k_1!k_2!(n-k_1-k_2)!}.\)
We have the following properties.\\
\((c1)\) There is a constant~\(C > 0\) such that for all~\(n \geq 3,\)~\(1 \leq l \leq N\) and~\(\frac{\eta_1 n}{2N} \leq k \leq \frac{2\eta_2 n}{N},\)
\begin{equation}\label{b_poi_diff}
|B(k;n,p_{l}) - Poi(k;np_{l})| \leq Poi(k;np_{l}) \left(1+ \frac{C n}{N^2}\right).
\end{equation}
\((c2)\) There is a constant~\(C > 0\) such that for all~\(n \geq 3,\)
and for any~\(1 \leq l_1,l_2 \leq N\) and~\(\frac{\eta_1 n}{2N} \leq k_1,k_2 \leq \frac{2\eta_2 n}{N},\)
\begin{eqnarray}
&&|B(k_1,k_2;n,p_{l_1},p_{l_2}) - Poi(k_1;np_{l_1})Poi(k_2; np_{l_2})| \nonumber\\
&&\;\;\;\;\leq Poi(k_1;np_{l_1}) Poi(k_2;np_{l_2})\left(1+ \frac{C n}{N^2}\right). \label{b_poi_diff2}
\end{eqnarray}
\emph{Proof of~\((c1)-(c2)\)}:
To prove~(\ref{b_poi_diff}) in~\((c1),\) we write~\(p_{l} =p\) for simplicity.
Use~\({n \choose k} \leq \frac{n^{k}}{k!}\) and~\(1-x \leq e^{-x}\) for~\(0 < x < 1\) to get
\[{n \choose k} p^{k} (1-p)^{n-k} \leq \frac{(np)^{k}}{k!}e^{-p(n-k)} = Poi(k;np)e^{kp}.\]
Using~(\ref{pl_def}) and the fact that~\(k \leq \frac{2\eta_2}{N}\) we get~\(e^{kp} \leq \exp\left(\frac{k\eta_2 n}{N}\right)
\leq \exp\left(2\eta_2 \frac{n}{N^2}\right)\)
and since~
\begin{equation}\label{ex_est_small_x}
e^{x} = 1+x+\sum \frac{x^{k}}{k!} \leq 1+x+ \sum_{k \geq 2}x^{k} \leq 1+2x
\end{equation}
for all~\(x\) small, we get~\(e^{kp} \leq 1 + \frac{4\eta_2 n}{N^2},\) proving the upper bound in~(\ref{b_poi_diff}).

To obtain a lower bound, we use the estimate
\begin{equation}\label{exp_low}
1-x \geq e^{-x-x^2}
\end{equation}
for all~\(0 < x< \frac{1}{2}.\) To prove~(\ref{exp_low}),
write \(\log(1-x) = -x - R(x)\) where~\[R(x) = \sum_{k \geq 2}\frac{x^{k}}{k} \leq \frac{1}{2}\sum_{k \geq 2}x^{k} =  \frac{x^2}{2(1-x)} \leq x^2\]
since~\(x < \frac{1}{2}.\)
Use~\({n \choose k} \geq \frac{(n-k)^{k}}{k!}\) and~(\ref{exp_low}) to get
\begin{equation}\label{bin2}
B(k;n,p) \geq \frac{1}{k!}(n-k)^{k} p^{k}e^{-p(n-k) - p^2(n-k)} = Poi(k;np)\left(1-\frac{k}{n}\right)^{k}e^{kp-(n-k)p^2}
\end{equation}
As before, using the fact that~\(\frac{\eta_1 n}{2N} \leq k \leq \frac{2\eta_2 n}{N}\) we get
\begin{equation}\label{temp_est2}
\left(1-\frac{k}{n}\right)^{k} \geq 1-\frac{k^2}{n} \geq 1-\frac{4\eta_2^2 n}{N^2}
\end{equation}
and~using~(\ref{pl_def}) we get
\begin{equation}
kp - (n-k)p^2 \geq kp - np^2 \geq \frac{\eta_1 n}{2N} \frac{\eta_1}{N}  - n\left(\frac{\eta_2 }{N}\right)^2 = -\eta \frac{n}{N^2}\label{temp_estt}
\end{equation}
where~\(\eta  = \eta^2_2 - \frac{\eta_1^2}{4}> 0,\) since~\(\epsilon_1 \leq  \epsilon_2\) and so~\(\eta_1  = \frac{\epsilon_1}{\epsilon_2} \leq \frac{\epsilon_2}{\epsilon_1} = \eta_2.\)
Using~(\ref{temp_est2}) and~(\ref{temp_estt}) into~(\ref{bin2}) gives
\begin{eqnarray}
B(k;n,p) &\geq& Poi(k;np)\left(1-\frac{\eta_1^2}{4} \frac{n}{N^2}\right) \exp\left(-\eta \frac{n}{N^2}\right) \nonumber\\
&\geq& Poi(k;np)\left(1-\frac{\eta_1^2}{4} \frac{n}{N^2}\right) \left(1-\eta \frac{n}{N^2}\right), \nonumber
\end{eqnarray}
since~\(e^{-x} \geq 1-x\) for~\(0 < x < 1.\) This proves~(\ref{b_poi_diff}).

To prove~(\ref{b_poi_diff2}), write~\(p_{l_1} =p_1, p_{l_2} = p_2\) and~\(B_{12} = B(k_1,k_2; n,p_1,p_2)\) for simplicity.
Use
\begin{equation}\label{eq_n_k1k2}
{n \choose k_1,k_2}  = \frac{1}{k_1!k_2!}n(n-1)\ldots (n-k_1-k_2+1) \leq \frac{n^{k_1+k_2}}{k_1!k_2!}
\end{equation}
to get
\begin{equation}\label{b12_est1}
B_{12} \leq \frac{(np_1)^{k_1}}{k_1!} \frac{(np_2)^{k_2}}{k_2!} e^{-(p_1+p_2)n}e^{(p_1+p_2)(k_1+k_2)}.
\end{equation}
Using~(\ref{pl_def}), we get~\(p_1+p_2 \leq \frac{2\eta_2 }{N}\) and since~\(k_1,k_2 \leq \frac{2\eta_2 n}{N}\)
we get using~(\ref{ex_est_small_x}) that
\begin{equation}\label{ep1p2}
e^{(p_1+p_2)(k_1+k_2)} \leq \exp\left(\frac{4\eta_2^2 n}{N^2}\right) \leq 1+\frac{8\eta_2^2 n}{N^2}
\end{equation}
for all~\(n\) large, since~\(\frac{n}{N^2} \longrightarrow 0\) as~\(n \rightarrow \infty\) (see~(\ref{N_est})).
Substituting~(\ref{ep1p2}) into~(\ref{b12_est1}), we get the upper bound for~\(B_{12}\) in~(\ref{b_poi_diff2}).

For the lower bound for~\(B_{12}\) again use~(\ref{eq_n_k1k2}) to get
\begin{equation}\nonumber
{n \choose k_1,k_2}  \geq \frac{1}{k_1!k_2!} (n-k_1-k_2)^{k_1+k_2} = \frac{n^{k_1+k_2}}{k_1!k_2!} \left(1-\frac{k_1+k_2}{n}\right)^{k_1+k_2}.
\end{equation}
Using~\((1-x)^{r} \geq 1-rx\) for~\(r,x > 0\) we further get
\begin{equation}\label{eq_n_k1k22}
{n \choose k_1,k_2} \geq \frac{n^{k_1+k_2}}{k_1!k_2!}  \left(1-\frac{(k_1+k_2)^2}{n}\right) \geq \frac{n^{k_1+k_2}}{k_1!k_2!}  \left(1-\frac{4\eta_2^2n}{N^2}\right)
\end{equation}
since~\(k_1,k_2 \leq \frac{2\eta_2n}{N}.\)
Substituting~(\ref{eq_n_k1k22}) into~(\ref{bin_dist2}) we get
\begin{equation}\label{b12_k12}
B_{12} \geq \frac{(np_1)^{k_1}}{k_1!} \frac{(np_2)^{k_2}}{k_2!} (1-p_1-p_2)^{n-k_1-k_2}\left(1-\frac{4\eta_2^2n}{N^2}\right).
\end{equation}

To evaluate~\((1-p_1-p_2)^{n-k_1-k_2},\) we use the estimate~(\ref{exp_low}) which is applicable since from~(\ref{pl_def}), we have~\(p_1+p_2 \leq \frac{2\eta_2}{N} \leq \frac{2\eta_2}{\sqrt{n}} \longrightarrow 0\) as~\(n \rightarrow \infty\) (see~(\ref{n_N})).
Using~(\ref{exp_low}), we get
\begin{equation}\label{p12_k12_est1}
(1-p_1-p_2)^{n-k_1-k_2} \geq e^{-(p_1+p_2)(n-k_1-k_2) - (p_1+p_2)^2(n-k_1-k_2)} = e^{-np_1} e^{-np_2} e^{I_{1}-I_2},
\end{equation}
where~\(I_{1}= (p_1+p_2)(k_1+k_2) \geq 0\) and~\(I_2 = (p_1+p_2)^2 (n-k_1-k_2) \leq  n(p_1+p_2)^2 \leq \frac{\eta_2^2n}{N^2}\)
for some constant~\(C_1 > 0,\)   by~(\ref{pl_def}). Using~\(e^{-x} \geq 1-x\) we get~\(e^{I_1-I_2} \geq e^{-I_2} \geq 1 - \frac{\eta_2^2 n}{N^2}\)
and so from~(\ref{p12_k12_est1}), we get
\begin{equation}\label{p12_k12_est2}
(1-p_1-p_2)^{n-k_1-k_2} \geq  e^{-np_1} e^{-np_2} \left(1-\frac{\eta_2^2 n}{N^2}\right).
\end{equation}
Using~(\ref{p12_k12_est2}) in~(\ref{b12_k12}), we get the lower bound for~\(B_{12}\) in~(\ref{b_poi_diff2}).~\(\qed\)

%WRTE MORE HERE... draw figures~\(\qed\)

Using properties~\((c1)-(c2)\) we prove Lemma~\ref{poi_bin_diff_lemm}.\\
\emph{Proof of~(\ref{tl_diff}) in Lemma~\ref{poi_bin_diff_lemm}}: Recall from~(\ref{nl_def}) that~\(N_l\) is the number of nodes
of the Binomial process~\(\{X_k\}\) in the square~\(S_l\) and let~\(U_l\) be the event as defined in~(\ref{ul_def}). Write~\(\mathbb{E}R_l = I_1 + I_2\)
where~
\begin{equation}\label{etl_1}
I_1 = \mathbb{E} R_l \ind(U_l) = \sum_{\frac{\eta_1 n}{2N} \leq k \leq \frac{2\eta_2 n}{N}} \mathbb{E} R_l \ind(N_l = k),I_2 = \mathbb{E}R_l \ind(U^c_l)
\end{equation}
and~\(\eta_1,\eta_2\) are as in~(\ref{pl_def}). Similarly~\(\mathbb{E}_0R^{(P)}_l = I^{(P)}_1 + I^{(P)}_2,\) where
\begin{equation}\label{et_poi1}
I^{(P)}_1 = \mathbb{E}_0(R^{(P)}_l \ind(U^{(P)}_l)),I^{(P)}_2 = \mathbb{E}_0(R^{(P)}_l \ind(U^{(P)}_l)^c),
\end{equation}
\(U^{(P)}_l = \left\{\frac{\eta_1 n}{2N} \leq N^{(P)}_l \leq \frac{2\eta_2 n}{N}\right\}\) is as defined in~(\ref{ul_def_p})
and~\(N^{(P)}_l\) is the number of nodes of the Poisson process~\({\cal P}\) inside the square~\(S_l\) (see discussion prior to~(\ref{poi_dist})).
From~(\ref{etl_1}) and~(\ref{et_poi1}), we therefore get
\begin{equation}\label{tl_diff33}
|\mathbb{E}R_l - \mathbb{E}_0R^{(P)}_l| \leq |I_1 - I^{(P)}_1| + I_2 + I^{(P)}_2.
\end{equation}

The remainder terms~\(I_2\) and~\(I^{(P)}_2\) satisfy
\begin{equation}\label{i2_est_fin2}
\max(I_2,I^{(P)}_2) \leq C (\mathbb{E} R_l) \frac{n}{N^2}
\end{equation}
for some constant~\(C > 0.\) We prove~(\ref{i2_est_fin2}) for~\(I_2\)
and an analogous proof holds for~\(I^{(P)}_2.\)
Indeed, every edge in the MST~\({\cal R}_l\)
containing all the nodes in the~\(r_n \times r_n\) square~\(S_l\)
has both endvertices within~\(S_l\) and so has length at most~\(r_n\sqrt{2}.\)
Since there are~\(N_l\) nodes in the square~\(S_l,\) there are~\(N_l-1 \leq N_l\)
edges in~\({\cal R}_l\) and so the length~\(R_l \leq N_l r_n\sqrt{2}\)
and
\begin{equation}\label{i2_etr22}
I_2  = \mathbb{E}R_l \ind(U_l^c) \leq r_n \sqrt{2} \mathbb{E}N_l\ind(U_l^c).
\end{equation}
Using the third expression in~(\ref{ul_2233}) to estimate~\(\mathbb{E}N_l \ind(U_l^c)\) we get
\begin{equation}\label{i2_etr2}
I_2  \leq C_1 r_n\sqrt{2} \frac{n}{N}\exp\left(-C_2 \frac{n}{N}\right) = C_1\sqrt{2} \left(r_n \sqrt{\frac{n}{N}}\right) \left(\sqrt{\frac{n}{N}}\exp\left(-C_2 \frac{n}{N}\right)\right)
\end{equation}
for some constants~\(C_1,C_2 > 0.\) From the lower bound in~(\ref{del_tn_b}) we have~\(\mathbb{E} R_l \geq C_3 r_n \sqrt{\frac{n}{N}}\)
and so
\begin{eqnarray}
I_2 \leq C_4 \left(\mathbb{E}R_l\right)\left(\sqrt{\frac{n}{N}}\exp\left(-C_2 \frac{n}{N}\right)\right) = C_4 \left(\mathbb{E}R_l\right)
\left(\frac{n}{N^2}\right) \delta_N \label{i2_etr33}
\end{eqnarray}
where
\begin{equation}\label{m_bds2}
\delta_N = \left(\frac{N^3}{n}\right)\exp\left(-\frac{C_2 n}{2N}\right) \leq \frac{n^2}{M^3(\log{n})^3} \exp\left(-\frac{C_2M}{2}\log{n}\right) \leq 1
\end{equation}
for all~\(n\) large, provided~\(M >0\) large. The first estimate in~(\ref{m_bds2}) follows from the upper bound~\(N \leq \frac{n}{M\log{n}}\) in~(\ref{n_N}). Fixing such an~\(M\) and using~(\ref{m_bds2}) in~(\ref{i2_etr33}), we get~(\ref{i2_est_fin2}).

%since~\(\frac{n}{N^2} \longrightarrow 0,\) we also have
%that~\(N \geq \sqrt{n}\) for all~\(n\) large. Thus
%\begin{equation}\label{m_bounds}
%\sqrt{n} \leq N \leq \frac{n}{M\log{n}}
%\end{equation}
%for all~\(n\) large and using~(\ref{m_bounds}) we have

To estimate the difference~\(I_1 - I^{(P)}_1\) in~(\ref{tl_diff33}),
recall that given~\(N_l = k,\) the nodes in~\(S_l\) are independently distributed in~\(S_l\) with distribution~\(\frac{f(.)}{\int_{S_l} f(x)dx}\) (see~(\ref{z_dist})) and so
\begin{equation}
I_1 = \sum_{\frac{\eta_1 n}{2N} \leq k \leq \frac{2\eta_2 n}{N}}\mathbb{P}(N_l = k)\mathbb{E}(R_l| N_l =k)  = \sum_{\frac{\eta_1 n}{2N} \leq k \leq \frac{2\eta_2 n}{N}} B(k; n,p_l) \Delta(k,q_l) \label{i1_est1}
\end{equation}
where~\(B(k;n,p_{l})\) is the Binomial probability distribution as defined in~(\ref{bin_dist}),~\(q_l = \int_{S_l} f(x) dx,\)
\begin{equation}\label{del_p_def}
\Delta(k,q_l) = \mathbb{E}(R_l|N_l = k) = \int_{S_l} MST(z_1,\ldots,z_k;S_l) \frac{f(z_1)}{q_l}\ldots \frac{f(z_k)}{q_l} dz_1\ldots dz_k
\end{equation}
and~\(MST(z_1,\ldots,z_k;S_l)\) is the length the MST containing all the nodes~\(z_1,\ldots,z_k \in S_l\) (see~(\ref{min_tree})).

Similarly, as argued in~(\ref{dlp_k}), given~\(N_l^{(P)} = k,\) the nodes of the Poisson process~\({\cal P}\) are
also distributed in~\(S_l\) according to distribution~\(\frac{f(.)}{\int_{S_l} f(x)dx}.\) Therefore~\(\mathbb{E}(R^{(P)}_l | N^{(P)}_l = k) = \Delta(k,q_l)\) as defined in~(\ref{del_p_def}) and so
\begin{equation}\label{ip_est11}
I^{(P)}_1 = \sum_{\frac{\eta_1 n}{2N} \leq k \leq \frac{2\eta_2 n}{N}}\Delta(k,q_l) Poi(k;np_l),
\end{equation}
where~\(Poi(k;np_l)\) is the Poisson distribution as defined in~(\ref{poi_dist}).
From~(\ref{i1_est1}) and~(\ref{ip_est11}), we therefore get
\begin{equation}\label{i1p_diff}
|I_1 - I^{(P)}_1| \leq \sum_{\frac{\eta_1 n}{2N} \leq k \leq \frac{2\eta_2 n}{N}} \Delta(k,q_l) |B(k;n,p_l) - Poi(k;np_l)|.
\end{equation}

Using estimate~(\ref{b_poi_diff}) of property~\((c1)\) to approximate the Binomial distribution with the Poisson distribution, we get
\begin{eqnarray}
|I_1-I^{(P)}_1| &\leq& C_1 \left(\sum_{\frac{\eta_1 n}{2N} \leq k \leq \frac{2\eta_2 n}{N}} Poi(k;np_l) \Delta(k,q_l) \right)\frac{n}{N^2} \nonumber\\
&\leq& C_1 \left(\sum_{k \geq 0}Poi(k;np_l) \Delta(k,q_l)\right) \frac{n}{N^2} \nonumber\\
&=& C_1 \left(\mathbb{E}_0(R^{(P)}_l)\right) \frac{n}{N^2} \label{i1_temp_estt}
\end{eqnarray}
for some constant~\(C_1 >0.\) But~\(\mathbb{E}_0(R^{(P)}_l)\) and~\(\mathbb{E}R_l\) both are bounded above and below by constant multiples of~\(r_n \sqrt{\frac{n}{N}}\) (see~(\ref{del_tn_b}) and~(\ref{del_tn})). From~(\ref{i1_temp_estt}), we therefore get
\begin{equation}\label{i1_fin}
|I_1-I^{(P)}_1| \leq C_3 \left(\mathbb{E} R_l\right) \frac{n}{N^2}
\end{equation}
for some constant~\(C_3 >0.\)
Substituting~(\ref{i1_fin}) and~(\ref{i2_est_fin2}) into~(\ref{tl_diff33}) gives
\[|\mathbb{E}R_l- \mathbb{E}_0R^{(P)}_l| \leq C_4 \left(\mathbb{E} R_l\right) \frac{n}{N^2} \leq C_5 \left(\frac{r_nn^{3/2}}{N^{5/2}}\right),\] for some positive constants~\(C_4,C_5,\) again using the upper bound for~\(\mathbb{E}R_l\) from~(\ref{del_tn_b}). This proves~(\ref{tl_diff}).~\(\qed\)

\emph{Proof of~(\ref{tl_diff2}) of Lemma~\ref{poi_bin_diff_lemm}}: The proof is analogous to~(\ref{tl_diff}).

Write~\(\mathbb{E} R_{l_1} R_{l_2} = J_1 + J_2\) where~\(J_1 = \mathbb{E} R_{l_1} R_{l_2} \ind( U_{l_1} \cap U_{l_2})\) and~\(J_2 = \mathbb{E}R_{l_1} R_{l_2} \ind(U^c_{l_1} \cup U^{c}_{l_2}).\)
Similarly, for the Poisson case let~\(U^{(P)}_l\) be the event defined in~(\ref{ul_def_p}) and define analogous terms~\(J_1^{(P)}\) and~\(J_2^{(P)}\) so that ~\(\mathbb{E}_0 R^{(P)}_{l_1} R^{(P)}_{l_2} = J^{(P)}_1 + J^{(P)}_2.\)
The difference
 \begin{equation}\label{r12_diff}
|\mathbb{E}R_{l_1}R_{l_2} - \mathbb{E}_0 R^{(P)}_{l_1} R^{(P)}_{l_2}| \leq |J_1 - J_1^{(P)}| + J_2 + J_2^{(P)}.
\end{equation}

Arguing as in~(\ref{i2_est_fin2}), the remainder terms~\(J_2\) and~\(J^{(P)}_2\) satisfy
\begin{equation}\label{i12_est_fin2}
\max(J_2,J^{(P)}_2) \leq C_1 (\mathbb{E} R_{l_1}\mathbb{E}R_{l_2}) \frac{n}{N^2} \leq C_2 \left(\frac{r_n^2 n^2}{N^3}\right)
\end{equation}
for some constants~\(C_1,C_2 > 0.\) We prove~(\ref{i12_est_fin2}) for~\(J_2\)
and an analogous proof holds for~\(J^{(P)}_2.\)
As argued in the proof of~(\ref{i2_est_fin2}), every
one of the~\(N_{l_1}\) edges in the MST~\({\cal R}_{l_1}\) of length~\(R_{l_1}\)
has both endvertices within~\(S_{l_1}\) and so has length at most~\(r_n\sqrt{2}.\)
Therefore
\begin{equation}\label{i12_etr22}
J_2  = \mathbb{E}R_{l_1}R_{l_2} \ind(U_{l_1}^c \cup U_{l_2}^c) \leq \left(r_n \sqrt{2}\right)^2 \mathbb{E}N_{l_1}N_{l_2}\ind(U^c_{l_1} \cup U^{c}_{l_2}).
\end{equation}
Using Cauchy-Schwarz inequality,
\begin{equation}\label{nl12_1}
\mathbb{E}N_{l_1}N_{l_2}\ind(U^c_{l_1} \cup U^{c}_{l_2}) \leq \left(\mathbb{E}N^2_{l_1}N^2_{l_2}\right)^{\frac{1}{2}}
\mathbb{P}\left(U^c_{l_1} \cup U_{l_2}^c\right)^{\frac{1}{2}} \leq \left(\mathbb{E}N^2_{l_1}N^2_{l_2}\right)^{\frac{1}{2}} \exp\left(-2C \frac{n}{N}\right)
\end{equation}
for some constant~\(C > 0\)  using the estimate~(\ref{ul_est}).

To evaluate~\(\mathbb{E} N^2_{l_1} N^2_{l_2},\) use~\(ab \leq \frac{a^2+b^2}{2}\) to write~\(\mathbb{E} N^2_{l_1} N^2_{l_2} \leq \frac{1}{2}\left(\mathbb{E}N_{l_1}^4 + \mathbb{E} N_{l_2}^4\right)\)
and use the fact that the term~\(N_l\) is Binomially distributed with parameters~\(n\) and~\(p_l,\) where~\(p_l \leq  \frac{\eta_2}{N}\) (see~(\ref{pl_def}))
and~\(\eta_2\) does not depend on~\(l\) or~\(n.\) Therefore
~\(\mathbb{E}N_l^4 \leq C_1 (np_l)^4 \leq C_2 \left(\frac{n}{N}\right)^4\) for some constants~\(C_1,C_2\) not depending on~\(l\) or~\(n\)
and so~\(\mathbb{E} N^2_{l_1} N^2_{l_2} \leq C_3 \left(\frac{n}{N}\right)^4.\) Therefore~\(\mathbb{E}N_{l_1}N_{l_2}\ind(U^c_{l_1} \cup U^{c}_{l_2}) \leq C_4\left(\frac{n}{N}\right)^2 e^{-2C\frac{n}{N}}\) (see~(\ref{nl12_1}))
and so from~(\ref{i12_etr22})
\begin{equation}\nonumber
J_2  \leq C_5 r_n^2 \left(\frac{n}{N}\right)^2 \exp\left(-2C\frac{n}{N}\right) = C_5 \left(\frac{r_n^2n^2}{N^3}\right) N\exp\left(-2C\frac{n}{N}\right).
\end{equation}
Since~\(N \leq \frac{n}{M\log{n}}\) (see~(\ref{n_N})) we have that~\(Ne^{-2C\frac{n}{N}} \leq \frac{n}{M\log{n}} e^{-2CM\log{n}} \leq 1\) for all~\(n\) large provided~\(M > 0\) is large. Fixing such an~\(M,\) we get~(\ref{i12_est_fin2}).

To evaluate the difference~\(J_1-J^{(P)}_1,\) recall from discussion prior to~(\ref{i1_est1}) that given~\(N_l = k,\)
the nodes of the Binomial process are distributed in the square~\(S_l\) with distribution~(\ref{z_dist}). Similarly,
given~\(N_l^{(P)} = k,\) the nodes of the Poisson process are also distributed according to~(\ref{z_dist}).
Therefore analogous to~(\ref{i1p_diff}) we get
\begin{equation}\label{j1_est}
|J_1-J^{(P)}_1| = \sum_{\frac{\eta_1 n}{2N} \leq k_1,k_2 \leq \frac{2\eta_2 n}{N}} |B_{l_1,l_2} - Poi(k_1;np_{l_1}) Poi(k_2;np_{l_2})| \Delta(k_1,q_{l_1}) \Delta(k_2,q_{l_2})
\end{equation}
where~\(q_{l_1}, q_{l_2}\) and~\(\Delta(.,.)\) are as defined in~(\ref{del_p_def}) and
\(B_{l_1,l_2} = B(k_1,k_2;n,p_{l_1},p_{l_2})\) is as defined in~(\ref{bin_dist2}).
Using~(\ref{b_poi_diff2}) and arguing as in~(\ref{i1_temp_estt}) we then get~\(|J_1-J^{(P)}_1| \leq C \mathbb{E}_0(R^{(P)}_{l_1}) \mathbb{E}_0(R^{(P)}_{l_2}) \left(\frac{n}{N^2}\right)\)
for some constant~\(C >0.\) Using the fact that bound~\(\mathbb{E}_0(R^{(P)}_{l_1})\) and~\(\mathbb{E}R_{l_1}\) are both bounded above and below by constant multiples of~\(r_n \sqrt{\frac{n}{N}}\) (see~(\ref{del_tn}) and~(\ref{del_tn_b})), we get~(\ref{tl_diff2}).~\(\qed\)

\emph{Proof of Lemma~\ref{cov_lemma}}:
Since Poisson process is independent on disjoint subsets,
we have \[cov_0(R^{(P)}_{l_1},R^{(P)}_{l_2}) = \mathbb{E}_0(R^{(P)}_{l_1}R^{(P)}_{l_2}) - \mathbb{E}_0R^{(P)}_{l_1}\mathbb{E}_0R^{(P)}_{l_2} = 0.\]

Therefore write
\[|cov(R_{l_1},R_{l_2})| = |cov(R_{l_1},R_{l_2}) - cov_0(R^{(P)}_{l_1},R^{(P)}_{l_2})| \leq Z_1 + Z_2 + Z_3,\]
where
\[Z_1 = |\mathbb{E} R_{l_1} R_{l_2} - \mathbb{E}_0 R^{(P)}_{l_1} R^{(P)}_{l_2}| \leq C\left(\frac{r_n^2 n^2}{N^3}\right),\]
\[Z_2 = |\mathbb{E}_0R^{(P)}_{l_1}\mathbb{E}_0R^{(P)}_{l_2} - \mathbb{E}R_{l_1}\mathbb{E}R_{l_2}| \leq Z_3 + Z_4,\]
\[Z_3 = |\mathbb{E}_0R^{(P)}_{l_1} - \mathbb{E}R_{l_1}|\mathbb{E}_0R^{(P)}_{l_2} \leq C\left(\frac{r_n n^{3/2}}{N^{5/2}}\right)
\left(r_n\sqrt{\frac{n}{N}}\right) = C \left(\frac{r_n^2 n^2}{N^3}\right)\]
and similarly,
\[Z_4 = \mathbb{E}R_{l_1}|\mathbb{E}_0R^{(P)}_{l_2} - \mathbb{E}R_{l_2}| \leq  C\frac{r_n^2 n^2}{N^3},\]
for some constant~\(C > 0.\) The estimate for~\(Z_1\) follows from~(\ref{tl_diff2}) and the estimates for~\(Z_3\) and~\(Z_4\)
follow from~(\ref{tl_diff}) and the estimates for~\(\mathbb{E}R_l\) and~\(\mathbb{E}_0R^{(P)}_l\) in~(\ref{del_tn_b}) and~(\ref{del_tn}),
respectively.~\(\qed\)

%WRTE MORE HERE... draw figures~\(\qed\)

\setcounter{equation}{0}
\renewcommand\theequation{\thesection.\arabic{equation}}
\section{Proof of Theorem~\ref{mst_thm}}\label{pf_mst}
For~\(1 \leq l \leq N,\) recall that~\(R_l\) is the length of the MST
containing all the~\(N_l\) nodes of~\(\{X_k\}\) present in the square~\(S_l.\)
The first step is to see that~\(MSTC_n\) is well approximated by~\(\sum_{l=1}^{N} R_l.\)
Recall that~\(s_n\) denotes the intercity distance i.e., the minimum distance
between squares in~\(\{S_l\}\)
(see paragraph prior to~(\ref{N_est})).

We have the following bounds for~\(MSTC_n.\)
\begin{Lemma}\label{mstc_lemm1} We have that
\begin{equation}\label{mstc_up}
MSTC_n \leq \left(V_n + (N-1)(s_n+8r_n)\right) \ind(U_{tot}(n)) + 3\sqrt{n} \ind(U^c_{tot}(n)),
\end{equation}
where
\begin{equation}\label{vn_def}
V_n := \sum_{l=1}^{N} R_l, U_{tot} = U_{tot}(n) := \bigcap_{l=1}^{N} U_l
\end{equation}
and~\(U_l\) is the event defined in~(\ref{ul_def}).
If the intercity distance~\(s_n > r_n \sqrt{2},\) then
\begin{equation}\label{mstc_low}
MSTC_n \geq V_n.
\end{equation}
\end{Lemma}

\emph{Proof of~(\ref{mstc_up}) of Lemma~\ref{mstc_lemm1}}: We construct a tree containing all the nodes\\
\(\{X_k\}_{1 \leq k \leq n}\)
and satisfying the upper bound in~(\ref{mstc_up}). Suppose that the event~\(U_{tot}\) occurs
so that each square~\(S_l, 1 \leq l \leq N\) contains at least
\begin{equation}\label{min_nodes}
\frac{\eta_1 n}{2N} \geq \frac{\eta_1 M}{2} \log{n} \geq 2,
\end{equation}
nodes of~\(\{X_k\}_{1 \leq k \leq n}\) for all~\(n\) large, by~(\ref{n_N}).
Let~\({\cal T}(l) \neq \emptyset\) be the MST containing all the nodes of~\(S_l.\)

Recall from the discussion following~(\ref{N_est}) that the cities are well connected in the sense that the vertices~\(\{z_l\}\)
corresponding to the centres of the squares~\(\{S_l\}\) is a connected
graph~\(G_Z \subset \mathbb{Z}^2.\)
The spanning tree~\(T_Z \subset G_Z\) contains~\(N-1\) edges~\(f_k, 1 \leq k \leq N-1.\)
Let~\(f_k\) have endvertices~\(z_1,z_2 \in \mathbb{Z}^2\)
and let~\(S_{l_1}\) and~\(S_{l_2}\) be the corresponding squares whose centres
are associated with~\(z_1\) and~\(z_2,\) respectively. Pick an edge~\(e_k\) with one endvertex
being a node of~\(\{X_k\}\) in~\(S_{l_1}\) and another endvertex being a node of~\(\{X_k\}\) in~\(S_{l_2}.\)
Performing this operation iteratively, we obtain~\(N-1\) edges~\(\{e_k\}_{1 \leq k \leq N-1}.\)

The union of the MSTs and the edges
\[{\cal T}_{up} := \bigcup_{l=1}^{N} {\cal T}(S_l) \bigcup \bigcup_{k=1}^{N-1} e_k\]
is a tree containing all the nodes~\(\{X_k\}_{1 \leq k\leq n}\)
and whose length is
\[L({\cal T}_{up}) = \sum_{l=1}^{N} R_l + \sum_{k=1}^{N-1} l(e_k) \leq \sum_{l=1}^{N} R_l + (N-1)(s_n+8r_n),\]
since each edge~\(e_k\) has length at most~\(s_n + 8r_n,\)
the sum of  the intercity distance and the total perimeter of the two~\(r_n \times r_n\)
squares containing the endvertices of~\(e_k.\)

%we obtain the upper bound in~(\ref{key_est_mst}) for the case when the event~\(U_{tot}\) occurs.

If the event~\(U_{tot}\) does not occur, then by the strips estimate~(\ref{mst_ab}),
the minimum spanning tree containing all the nodes~\(\{X_k\}_{1 \leq k \leq n}\) has length
at most~\(3\sqrt{n}.\)\(\qed\)

To prove the lower bound~(\ref{mstc_low}) in Lemma~\ref{mstc_lemm1}, we need additional properties.
Recall from~(\ref{len_cyc_def}) that~\({\cal T}_n\)
is the minimum spanning tree containing all the nodes~\(\{X_k\}_{1 \leq k \leq n}.\)
Suppose there are two nodes~\(a,b \in \{X_k\}\) in some square~\(S_l, 1 \leq l \leq N.\)
Since~\({\cal T}_n\) is a tree, there is a unique path~\({\cal P}_{ab} \subset {\cal T}_n\) containing~\(a\) and~\(b\) as endvertices.
The following crucial property also holds.\\
\((g1)\) Every node in~\({\cal P}_{ab}\) belongs to the square~\(S_l.\)\\
\emph{Proof of~\((g1)\)}: We prove by contradiction and
suppose that the path~\({\cal P}_{ab}\) contains a node outside the square~\(S_l.\)
This means that~\({\cal P}_{ab}\) ``exits" and ``re-enters" the square~\(S_l\) at two distinct nodes.
Without loss of generality, we assume that~\(a\) and~\(b\) are the exit and entry points; i.e.,
there are edges~\(e_a\) and~\(e_b\) both in~\({\cal P}_{ab}\) such that~\(e_a\) contains~\(a\)
as an endvertex  and~\(e_b\) contains~\(b\) as an endvertex.

\begin{figure}[tbp]
\centering
%\fbox{
\includegraphics[width=2in, trim= 60 200 150 150, clip=true]{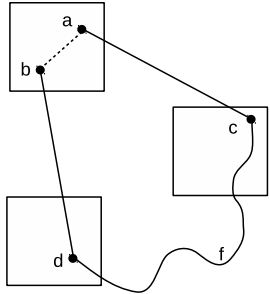}
%}
\caption{Modifying the path~\(P_{ab} = acfdb\) to obtain a new tree~\({\cal T}_{new}.\)}
\label{mst_fig1}
\end{figure}

If~\(c\) and~\(d\)
are the other endvertices of~\(e_a\) and~\(e_b\) respectively, then~\(c\) and~\(d\) both lie
outside~\(S_l,\) as shown in Figure~\ref{mst_fig1}. Here, the path~\({\cal P}_{ab}=acfdb\) is the union
of the two edges~\(ac, bd\) and the wavy path~\(cfd.\)

Since the distance between any two squares in~\(\{S_j\}\)
is at least~\(s_n > r_n\sqrt{2} ,\) the edges~\(ac\) and~\(bd\)
have length at least~\(s_n > r_n\sqrt{2},\) each. The edge~\(ab\) however has length
at most~\(r_n \sqrt{2}.\) Consider the new graph~\({\cal T}_{new}\) formed by deleting
the edge~\(ac\) and adding the edge~\(ab.\)
The graph~\(T_{new}\) is a tree and by construction, the sum of the length of edges
in~\({\cal T}_{new}\) is strictly less than
the sum of length of edges in the MST~\({\cal T}_n.\) This is a contradiction and so
all nodes of~\({\cal P}_{ab}\) are contained in the square~\(S_l.\)~\(\qed\)

\emph{Proof of~(\ref{mstc_low}) in Lemma~\ref{mstc_lemm1}}: For~\(1 \leq l \leq N,\)
let~\({\cal H}_n(l)\) be the subgraph of~\({\cal T}_n\)
containing all the nodes of~\(S_l\) and all edges with both endvertices inside~\(S_l.\)
From property~\((g1),\) the graph~\({\cal T}_n(l)\) is connected and is therefore a tree.
The length of~\({\cal T}_n(l)\) is at least~\(R_l,\)
the length of the MST containing all the nodes of~\(S_l.\)
Since the above statement is true for each~\(1 \leq l \leq N,\) we obtain the lower bound in~(\ref{mstc_low}).~\(\qed\)

We use Lemma~\ref{mstc_lemm1} to prove Theorem~\ref{mst_thm}.
From Lemma~\ref{mstc_lemm1},
we have that the overall minimum length~\(MSTC_n\) is bounded above and below by the sum of
the local MST lengths~\(\sum_{l=1}^{N}R_l\) apart from some residual terms.
From the bounds on~\(\mathbb{E}R_l\) in~(\ref{del_tn}) of Lemma~\ref{rl_lemma},
we have that~\(\sum_{l=1}^{N} \mathbb{E}R_l\) is of the order of~\(N r_n \sqrt{\frac{n}{N}}  = r_n \sqrt{nN} = b_n\)
as defined in~(\ref{bn_def}). We therefore study the convergence of~\(\frac{MSTC_n}{b_n}.\)
We henceforth fix~\(M > 0\) large so that~(\ref{cov_tl_est}) of Lemma~\ref{cov_lemma} holds.\\\\
\emph{Proof of~(\ref{conv_mst_prob}) in Theorem~\ref{mst_thm}}:
From the upper and lower bounds~(\ref{mstc_up}) and~(\ref{mstc_low}) in Lemma~\ref{mstc_lemm1},
we have that
\begin{eqnarray}
\frac{1}{b_n}(V_n - \mathbb{E}V_n)  - \Delta_n \leq \frac{1}{b_n}\left(MSTC_n - \mathbb{E} MSTC_n\right) \leq \frac{1}{b_n}(V_n - \mathbb{E}V_n) + \Delta_n \label{key_est2}
\end{eqnarray}
where~\(V_n = \sum_{l=1}^{N} R_l\) is as defined~(\ref{vn_def}) and
\begin{equation}\label{delta_n_def}
\Delta_n = \frac{2(N-1)(s_n+8r_n)}{b_n}\ind(U_{tot}(n)) + \frac{4\sqrt{n}}{b_n}\ind(U^c_{tot}(n)) \nonumber.
\end{equation}

The variance of~\(V_n\) satisfies
\begin{equation}\label{var_vn}
var(V_n) \leq C \frac{r_n^2 n^2}{N} = C b_n^2 \left(\frac{n}{N^2}\right)
%var(V_n) \leq C (\log{n})^2
\end{equation}
for some constant~\(C > 0\) and all~\(n\) large and  since~\(\frac{n}{N^2} \longrightarrow 0\) (see~(\ref{N_est})),
we get that
\begin{equation}\label{vn_prob_conv}
\frac{1}{b_n}\left(V_n - \mathbb{E}V_n\right) \longrightarrow 0 \text{ in probability }
\end{equation}
as~\(n \rightarrow \infty.\) Also
\begin{equation}\label{del_n_conv}
\Delta_n \longrightarrow 0 \text{ a.s.}
\end{equation}
as~\(n \rightarrow \infty.\) This proves~(\ref{conv_mst_prob}) and we prove~(\ref{var_vn}) and~(\ref{del_n_conv}) separately below.

\emph{Proof of~(\ref{var_vn})}: Write
\begin{eqnarray}
var(V_n) &=& \sum_{l} var(R_l) + \sum_{l_1,l_2} cov(R_{l_1},R_{l_2}) \nonumber\\
&\leq& \sum_{l} \mathbb{E}R_l^2  + \sum_{l_1,l_2} cov(R_{l_1},R_{l_2}), \label{var_vn_est1}
\end{eqnarray}
where~\(cov(X,Y) = \mathbb{E}XY - \mathbb{E}X\mathbb{E}Y.\)  Using~(\ref{del_tn}) of Lemma~\ref{rl_lemma} to
estimate~\(\mathbb{E}R_l^2\) we get~
\begin{equation}\label{f_t}
\sum_{l=1}^{N} \mathbb{E} R_l^2 \leq N C_1 \left(r_n \sqrt{\frac{n}{N}}\right)^2 = C_1r_n^2 n
\end{equation}
for some constant~\(C_1 > 0.\) Similarly using estimate~(\ref{cov_tl_est}) of Lemma~\ref{cov_lemma} for the covariance,
we get
\begin{equation}\label{s_t}
\sum_{l_1,l_2} cov(R_{l_1},R_{l_2})\leq N^2 \left(C_2 \frac{r_n^2 n^2}{N^3}\right) = C_2 \frac{r_n^2 n^2}{N}.
\end{equation}
for some constants~\(C > 0.\) Substituting~(\ref{f_t}) and~(\ref{s_t}) into~(\ref{var_vn_est1}), we get
\[var(V_n) \leq C_1 r_n^2 n + C_2 \frac{r_n^2 n^2}{N} = \frac{r_n^2 n^2}{N}\left(C_1 \frac{N}{n} + C_2\right).\]
Since~\(\frac{N}{n} \leq \frac{1}{M\log{n}} \leq 1\) for all~\(n\) large~(see~(\ref{n_N})), we get that~\(var(V_n) \leq C_3\frac{r_n^2 n^2}{N}\)
for some positive constant~\(C_3\) and for all~\(n\) large.

\emph{Proof of~(\ref{del_n_conv})}: From~(\ref{delta_n_def}) and the fact that~\(r_n < r_n\sqrt{2} < s_n\) (see statement of the Theorem),
we get~
\begin{equation}\label{delta_n}
0 \leq \Delta_n \leq \frac{18N s_n}{b_n}  + \frac{4\sqrt{n}}{b_n} \ind(U^c_{tot}(n))
\end{equation}
and so
\begin{equation}\label{delta_n2}
0 \leq \limsup_n \Delta_n \leq \limsup_n \frac{4\sqrt{n}}{b_n} \ind(U^c_{tot}(n)),
\end{equation}
since \(\frac{Ns_n}{b_n} \longrightarrow 0\) as~\(n \rightarrow \infty\) by the statement of the Theorem.
From the estimate for the event~\(U_l\) in~(\ref{ul_def}),
\begin{equation}\label{u_tot_est}
\mathbb{P}(U^c_{tot}(n)) \leq \sum_{l=1}^{N} \mathbb{P}(U^c_l) \leq N \exp\left(-C\frac{n}{N}\right),
\end{equation}
for some constant~\(C > 0.\) Using the fact that~\(\frac{n}{N} \geq M\log{n}\) (see~(\ref{n_N})), we get
\begin{equation}\label{u_tot_est2}
\mathbb{P}(U^c_{tot}(n)) \leq \frac{n}{M\log{n}} \frac{1}{n^{MC}} \leq \frac{1}{n^2},
\end{equation}
provided~\(M > 0\) is large. Fixing such an~\(M,\) we have from Borell-Cantelli lemma that~\(\mathbb{P}(\limsup_n U^{c}_{tot}(n))=0\) and so a.s. \(\ind(U_{tot}^c(n)) = 0\) for all large~\(n.\) From~(\ref{delta_n2}), we therefore get~(\ref{del_n_conv}).~\(\qed\)

\emph{Proof of~(\ref{exp_mstc}) in Theorem~\ref{mst_thm}}: Recalling that~\(V_n = \sum_{i=1}^{N} R_l\)
from~(\ref{vn_def}), we use Lemma~\ref{mstc_lemm1}
to get
\begin{equation}\label{vn_est_tl}
\mathbb{E}V_n \leq \mathbb{E} MSTC_n \leq \mathbb{E}V_n + b_n\mathbb{E}\Delta_n,
\end{equation}
where~\(\Delta_n\) satisfies (see~(\ref{delta_n}))
\begin{equation}\label{edel_n}
\mathbb{E}\Delta_n \leq \frac{18Ns_n}{b_n}  + \frac{4\sqrt{n}}{b_n} \mathbb{P}(U_{tot}^c(n)) \leq 18 + \frac{4\sqrt{n}}{b_n} \mathbb{P}(U_{tot}^c(n)),
\end{equation}
since~\(\frac{Ns_n}{b_n} \longrightarrow 0\) as~\(n \rightarrow \infty\) (see statement of the Theorem).
Using~(\ref{u_tot_est2}) for estimating the probability of the event~\(U_{tot}\) we get
\begin{equation}\label{u_tot_bn}
\sqrt{n}\mathbb{P}(U^{c}_{tot}(n)) \leq \frac{\sqrt{n}}{n^2} \leq \sqrt{\frac{M\log{n}}{n}} \leq r_n \leq r_n \sqrt{nN}  = b_n
\end{equation}
for all~\(n\) large, where the second inequality is true by the condition for~\(r_n\) in~(\ref{N_est}). %provided~\(M >0\) is large. Fix such an~\(M.\)
Thus~\(\frac{4\sqrt{n}}{b_n}\mathbb{P}(U^{c}_{tot}(n))  \leq 4\)
and so~\(\mathbb{E}\Delta_n \leq 22\) and
\begin{equation}\label{vn_est_tl2}
\mathbb{E}V_n \leq \mathbb{E} TSPC_n \leq \mathbb{E}V_n + 22b_n,
\end{equation}
by~(\ref{edel_n}) and~(\ref{vn_est_tl}), respectively.

To estimate~\(\mathbb{E}V_n\) use the bounds for~\(\mathbb{E}R_l\) in~(\ref{del_tn}) of Lemma~\ref{rl_lemma} to get
\begin{equation}\label{vn_estb4}
C_1 b_n = N\left(C_1 r_n \sqrt{\frac{n}{N}}\right) \leq \mathbb{E}V_n \leq N\left(C_2 r_n \sqrt{\frac{n}{N}}\right) = C_2 b_n
\end{equation}
for some constants~\(C_1,C_2  >0.\) From~(\ref{vn_estb4}) and~(\ref{vn_est_tl2}), we get the bounds for~\(\mathbb{E}MSTC_n\) in~(\ref{exp_mstc}).~\(\qed\)

\emph{Proof of~(\ref{eq_mst2}) of Theorem~\ref{mst_thm}}:
We consider Poissonization and recall the Poisson process~\({\cal P}\) on the squares~\(\{S_l\}_{1 \leq l \leq N},\)
defined on the probability space~\((\Omega_0, {\cal F}_0, \mathbb{P}_0)\) (see paragraph prior to~(\ref{poi_dist})).
Analogous to~\(MSTC_n\) defined in~(\ref{min_weight_tree}), let~\(MSTC^{(P)}_n\) denote the length of the MST containing all
the nodes of the Poisson process~\({\cal P}.\) Recall from~(\ref{rl_def_p}) that~\(R^{(P)}_l\) denotes the length
of the MST containing all the nodes of~\({\cal P}\) in the square~\(S_l.\)

Analogous to~(\ref{mstc_low}), we have that
if the intercity distance~\(s_n > r_n \sqrt{2},\) then
\begin{equation}\label{tspc_low_p}
MSTC^{(P)}_n \geq V^{(P)}_n = \sum_{l=1}^{N} R_l^{(P)}.
\end{equation}

Define the event~\(E^{(P)}_l = \left\{R^{(P)}_l \geq \delta_4 r_n \sqrt{\frac{n}{N}}\right\},\) where~\(\delta_4\) is the constant
in~(\ref{tn_prob_est}) of Lemma~\ref{rl_lemma_poiss}. Since the Poisson process
is independent on disjoint sets, the events~\(E^{(P)}_l\) are independent
and each occurs with probability at least~\(\delta_5,\) by~(\ref{tn_prob_est}).
If~\(F^{(P)}_{sum} := \sum_{l=1}^{N} \ind(E^{(P)}_l)\) then~\(\mathbb{E}_0 \left(F^{(P)}_{sum}\right) \geq \delta_5 N\)
and from the standard Chernoff bound estimate for sums of independent Bernoulli
random variables (see Corollary \(A.1.14,\) pp. 312 of Alon and Spencer~(2008))
we also have~\(\mathbb{P}_0\left(F^{(P)}_{sum} \geq C_1 N\right)
\geq 1- e^{-2C_2 N}\)
for some positive constants~\(C_1\) and~\(C_2.\) If~\(F^{(P)}_{sum} \geq C_1N,\)
then \(\sum_{l=1}^{N} R^{(P)}_l \geq C_1 N \left(\delta_4 r_n \sqrt{\frac{n}{N}}\right) = C_3 b_n\)
for some constant~\(C_3 > 0\) and so from~(\ref{tspc_low_p}),
\begin{equation}\label{poi_t_up}
\mathbb{P}_0(MSTC^{(P)}_n \geq C_3 b_n) \geq 1- e^{-2C_2 N}
\end{equation}
for all~\(n\) large.

To convert the probability estimates to the Binomial process, let~\[A_P = \{TSPC^{(P)}_n \geq C_3 b_n\}, A = \{TSPC_n \geq C_3 b_n\}\]
and use the dePoissonization formula
\begin{equation}\label{de_poiss_ax}
\mathbb{P}(A) \geq 1- D \sqrt{n} \mathbb{P}(A^c_P)
\end{equation}
for some constant~\(D > 0\) and~(\ref{poi_t_up}) to get that
\begin{equation}\nonumber
\mathbb{P}(MSTC_n \geq C_3 b_n) \geq 1- D\sqrt{n} e^{-2C_2 N} = 1-e^{-  \alpha_N},
\end{equation}
where~\(\alpha_N = 2C_2N - \log{D} - \frac{1}{2}\log{n} \geq C_2N\) for all~\(n\) large, since~\(N \geq \sqrt{n}\) for all~\(n\) large~(see~(\ref{n_N})).
This proves~(\ref{eq_mst2}) and it only remains to prove~(\ref{de_poiss_ax}).

To prove~(\ref{de_poiss_ax}), let~\(N_P\) denote the random number of nodes of~\({\cal P}\) in all the squares~\(\cup_{j=1}^{N}S_j\) so
that~\(\mathbb{E}_0 N_P = n\) and~\(\mathbb{P}_0(N_P=n) = e^{-n}\frac{n^{n}}{n!} \geq \frac{D_1}{\sqrt{n}}\) for
some constant~\(D_1 > 0,\) using the Stirling formula. Given~\(N_P = n,\) the nodes of~\({\cal P}\)
are i.i.d.\ with distribution~\(g_N\) as defined in~(\ref{gn_def}); i.e., \(\mathbb{P}_0(A_P^c|N_P = n)  = \mathbb{P}(A^c)\) and so
\[\mathbb{P}_0(A_P^c) \geq \mathbb{P}_0(A_P^c|N_P = n) \mathbb{P}_0(N_P = n) =
\mathbb{P}(A^c) \mathbb{P}_0(N_P = n) \geq \mathbb{P}(A^c)\frac{D_1}{\sqrt{n}},\] proving~(\ref{de_poiss_ax}).~\(\qed\)

\emph{Proof of~(\ref{eq_mst3}) of Theorem~\ref{mst_thm}}:
As in the proof of~(\ref{eq_mst2}) above, we consider the Poisson process~\({\cal P}\)
on the squares~\(\{S_l\}_{1 \leq l \leq N}\) defined in the paragraph prior to~(\ref{poi_dist}).
As before, let~\(MSTC^{(P)}_n\) denote the length of the minimum length cycle containing all
the nodes of the Poisson process~\({\cal P}.\) Recall from~(\ref{rl_def_p}) that~\(R^{(P)}_l\) denotes the length
of the minimum length cycle containing all the nodes of~\({\cal P}\) in the square~\(S_l.\)

Analogous to~(\ref{mstc_up}) of Lemma~\ref{mstc_lemm1}, we have
\begin{equation}\label{tspc_up_p}
MSTC^{(P)}_n \leq \left(V^{(P)}_n + (N-1)(s_n+8r_n)\right) \ind(U^{(P)}_{tot}(n)) + 4\sqrt{n} \ind(U^{(P)}_{tot}(n))^c,
\end{equation}
where
\begin{equation}\label{vn_def1_p}
V^{(P)}_n := \sum_{l=1}^{N} R^{(P)}_l, U^{(P)}_{tot} = U^{(P)}_{tot}(n) := \bigcap_{l=1}^{N} U^{(P)}_l
\end{equation}
and~\(U^{(P)}_l = \{\frac{\eta_1 n}{2N} \leq N^{(P)}_l \leq \frac{2\eta_2 n}{N}\}\) is the event defined in~(\ref{ul_def_p}).
Recall that~\(N^{(P)}_l\) is the total number of nodes of~\({\cal P}\) inside the square~\(S_l.\)

Suppose now that the event~\(U_{tot}^{(P)}(n)\) occurs so that
\begin{equation}\label{tspc_ph}
MSTC_n^{(P)} \leq V_n^{(P)}  + (N-1)(s_n + 8r_n) = \sum_{l=1}^{N} R^{(P)}_l + (N-1)(s_n+8r_n).
\end{equation}
Since~\(U^{(P)}_l \supseteq U_{tot}^{(P)}\) occurs for every~\(1 \leq l \leq N,\)
we use the strips estimate~(\ref{mst_ab}) with~\(a = \frac{2\eta_2 n}{N}\) and~\(b = r_n\)
to get that the corresponding minimum length~\(R_l^{(P)} \leq 4b\sqrt{a} \leq C r_n \sqrt{\frac{n}{N}}\)
for some constant~\(C > 0\) and for every~\(1 \leq l \leq N.\)
Thus~\(V^{(P)}_n = \left(\sum_{l=1}^{N} R^{(P)}_l\right) \leq C b_n\)
and from~(\ref{tspc_ph})
we therefore get
\begin{equation}\label{tspc_ph1}
MSTC_n^{(P)} \leq C b_n + 2(N-1)(s_n + 8r_n) \leq Cb_n + 18Ns_n \leq (C+1)b_n,
\end{equation}
for all~\(n\) large. The second inequality in~(\ref{tspc_ph1}) is true since~\(r_n < r_n \sqrt{2} < s_n.\)
The final inequality in~(\ref{tspc_ph1}) is true since~\(\frac{Ns_n}{b_n} \longrightarrow 0\) and
so~\(\frac{Ns_n}{b_n} \leq \frac{1}{18}\) for all~\(n\) large.

Summarizing, we have that if the event~\(U_{tot}^{(P)}\) occurs,
then the overall minimum length~\(MSTC_n^{(P)} \leq C_1 b_n\)
for some constant~\(C_1 > 0.\)
To evaluate~\(\mathbb{P}(U^{(P)}_{tot}),\) use the estimate~(\ref{ul_est_p}) for the event~\(U^{(P)}_l\) to get
\begin{equation}\label{u_tot_p_est}
\mathbb{P}_0(U^{(P)}_{tot}) \geq 1- N\exp\left(-2C\frac{n}{N}\right)
\end{equation}
for some constant~\(C > 0.\) Thus
\begin{equation}\label{tutu}
\mathbb{P}_0\left(MSTC_n^{(P)} \leq C_1 b_n \right) \geq \mathbb{P}(U^{(P)}_{tot}) \geq 1-N\exp\left(-2C\frac{n}{N}\right).
\end{equation}

To convert the probabilities to the Binomial process, we again use
the dePoissonization formula~(\ref{de_poiss_ax})
to get that
\begin{equation}\label{tutu_bin}
\mathbb{P}\left(MSTC_n \leq C_1 b_n \right) \geq 1-DN\sqrt{n}\exp\left(-2C\frac{n}{N}\right) = 1-e^{-\delta_N},
\end{equation}
where~\(D > 0\) is as in~(\ref{de_poiss_ax}) and~\(\delta_N = 2C \frac{n}{N} - \log{D}-\log{N} - \frac{1}{2}\log{n}.\)
Since~\(\frac{n}{N} \geq M\log{n}\) for all~\(n\) large (see~(\ref{n_N})), we get
\[\log{D} + \log{N} + \frac{1}{2}\log{n} \leq \log{D} + \log\left(\frac{n}{M\log{n}}\right) + \frac{1}{2}\log{n} \leq 2\log{n} \leq C\frac{n}{N},\]
provided~\(M >0\) is large. Fixing such an~\(M\) we get that~\(\delta_{N} \geq C\frac{n}{N}\) and so~(\ref{eq_mst3}) follows
from~(\ref{tutu_bin}).~\(\qed\)

\setcounter{equation}{0}
\renewcommand\theequation{\thesection.\arabic{equation}}
\section{Proof of Theorem~\ref{var_mst_thm}}\label{pf_mst_un}
To prove Theorem~\ref{var_mst_thm}, we need a preliminary estimate regarding the
difference in the total length of the MSTs upon adding or deleting a single node.
For~\(n \geq 1,\) divide the unit square~\(S\) into~\(r_n \times r_n\) squares~\(\{S_i\}_{1 \leq i \leq N}\) each of side length~\(r_n\)
satisfying
\begin{equation}\label{rn_def_mst}
\frac{2M\log{n}}{n} \leq r_n^2 := \frac{2M\log{n} + c_n}{n} \leq \frac{3M\log{n}}{n},
\end{equation}
where~\(M\) is a large integer to be determined later
and~\(c_n \in (0,1)\) is chosen
such that~\(\frac{1}{41r_n}\) is an integer.

For~\(1 \leq i \leq N,\) let~\(N_i\) be the random number of nodes of~\(\{X_k\}_{1 \leq k \leq n}\)
in the square~\(S_i.\) Using~(\ref{f_eq}), the average number of nodes
\[\mathbb{E}N_i = n\int_{S_i} f(x) dx\]
satisfies
\begin{equation} \label{ave_per_sq}
8 \leq 2\epsilon_1 M \log{n} \leq n\epsilon_1 r_n^2 \leq \mathbb{E} N_i \leq n\epsilon_2 r_n^2  \leq 3\epsilon_2 M \log{n}
\end{equation}
where~\(\epsilon_1,\epsilon_2 > 0\) is as in~(\ref{f_eq}). The first estimate in~(\ref{ave_per_sq}) is true provided the constant~\(M > 0\)
is large and we fix such an~\(M\) henceforth. The other estimates in~(\ref{ave_per_sq}) follow from~(\ref{rn_def_mst}).

For~\(1 \leq j \leq n+1\) and~\(1 \leq i \leq N,\) let~\(Z_j(i)\) be the event
that the square~\(S_i\) contains between~\(\epsilon_1 M \log{n}\) and~\(4\epsilon_2 M\log{n}\)
nodes of~\(\{X_k\}_{1 \leq k \neq j \leq n+1}\) and define
\begin{equation}\label{z_tot_def}
Z_{tot}(n+1) := \bigcap_{1 \leq j \leq n+1} \bigcap_{i=1}^{N} Z_j(i).
\end{equation}
By standard Binomial estimates and~(\ref{ave_per_sq}) (see Corollary A.1.14, pp. 312, Alon and Spencer (2008))
\begin{equation}\label{zji_est}
\mathbb{P}(Z_j(i)) \geq 1-e^{-2C_1 M\log{n}}
\end{equation}
for some positive constant~\(C_1\) not depending no~\(i\) or~\(j.\) Thus
\[\mathbb{P}(Z_{tot}(n+1)) \geq 1- (n+1)\cdot N \cdot e^{-2C_1 M\log{n}}\]
and since the number of squares is~\(N = \frac{1}{r_n^2} \leq \frac{C_2 n}{\log{n}}\) for
some constant~\(C_2 > 0\) (see~(\ref{rn_def_mst})),
we get
\begin{equation}\label{z_tot_est}
\mathbb{P}(Z_{tot}(n+1)) \geq 1- \frac{C_2 n}{\log{n}}(n+1)e^{-2C_1 M\log{n}} \geq 1-e^{-C_1 M\log{n}}
\end{equation}
for all~\(n\) large, provided~\(M > 0\) is large. Fix such a~\(M.\)

Recall that~\({\cal T}_n\) is the MST containing all the~\(n\) nodes~\(\{X_k\}_{1 \leq k \leq n}.\)
The following Lemma estimates the edge lengths in the MSTs~\({\cal T}_{n+1}\) and~\({\cal T}_n.\)
\begin{Lemma}\label{one_point_est}
For~\(1 \leq j \leq n+1,\) let~\(MST_n(j)\) be the length of the minimal spanning tree
containing the nodes~\(\{X_k\}_{1 \leq k \neq j \leq n+1}.\) The difference
\begin{equation}\label{mst_diff}
|MST_{n+1} - MST_n(j)| \leq C_1 r_n \log{n} \ind(Z_{tot}(n+1)) + n\sqrt{2} \ind(Z^c_{tot}(n+1)),
\end{equation}
for some constant~\(C_1 > 0\) not depending on~\(j.\) Also, if~\(M > 0\) is large then
\begin{equation}\label{mst_diff_one}
\mathbb{E}|MST_{n+1} - MST_n| \leq C_2\frac{(\log{n})^{3/2}}{\sqrt{n}}
\end{equation}
for some constant~\(C_2 > 0.\)
\end{Lemma}
We henceforth fix~\(M\) large enough so that~(\ref{mst_diff_one}) is also satisfied.

We first perform some preliminary computations. For a square~\(S_i,\)
let~\({\cal N}_1(S_i)\) be the set of all squares in~\(\{S_l\}\) sharing a corner with~\(S_i.\)
For~\(k \geq 2,\) let~\({\cal N}_k(S_i)\) be the set of squares sharing a corner with
some square in~\({\cal N}_{k-1}(S_i).\) We use the following property to
prove Lemma~\ref{one_point_est}.\\
\((h1)\) Suppose the event~\(Z_{tot}(n+1)\) occurs and
suppose~\(X_j  = v \in S_i\) for some~\(1 \leq i \leq N.\)
Let~\(e\) be any edge in the tree~\({\cal T}_{n+1}\)
containing~\(v\) as an endvertex. If~\(u\) denotes the other
endvertex of~\(e,\) then~\(u \in S_k\) for some~\(S_k \in {\cal N}_{20}(S_i)\)
and the length of~\(e\) is at most~\(20r_n \sqrt{2}.\) \\
\emph{Proof of~\((h1)\)}: The fact that the edge length is at most~\(20r_n \sqrt{2}\)
is a consequence of the definition of~\({\cal N}_{20}(S_i).\)

We prove by contradiction and assume that~\(u\) does not lie in any square of~\({\cal N}_{20}(S_i).\)
Let~\(S_k \in {\cal N}_{10}(S_i)\) be a square
whose centre is at a distance of at least~\(5r_n\) from the centre
of~\(S_i,\) intersecting the edge~\((u, v).\) %as shown in Figure~\ref{mst_fig2}.
Since the event~\(Z_{tot}(n+1)\) occurs, the square~\(S_k\) contains a vertex~\(z\) which
also belongs to the MST~\({\cal T}_{n+1}.\) The distance between~\(z\) and~\(u\) is
strictly less than the distance between~\(v\) and~\(u.\)
Similarly the distance between~\(v\) and~\(z\) is strictly less than
the distance between~\(v\) and~\(u.\)

\begin{figure}
\centering
\begin{subfigure}{0.5\textwidth}
\centering
\includegraphics[width=3in, trim= 20 380 50 80, clip=true]{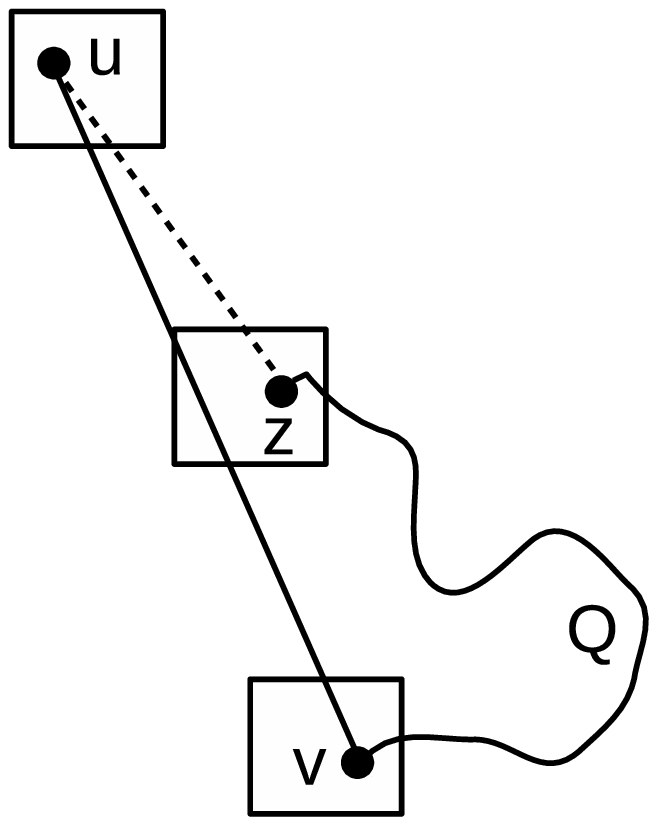}
   \caption{The node~\(u \notin {\cal P}_{vz} = vQz.\)}%\nonumber%{t1}
\end{subfigure}% seems this is important....
\begin{subfigure}{0.5\textwidth}
\centering
   \includegraphics[width=3in, trim= 20 380 50 80, clip=true]{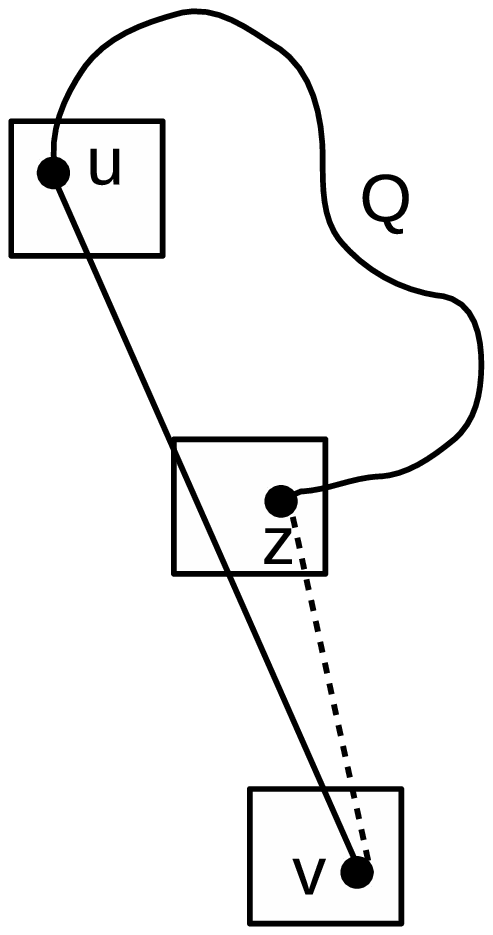}
   \caption{The node~\(u \in {\cal P}_{vz} = vuQz.\)} %\nonumber%{t2}
  \end{subfigure}

\caption{Creating the new tree~\({\cal T}_{new}\) depending on whether~\(u \in {\cal P}_{vz}\) or not.}
\label{mst_fig2}
\end{figure}

Let~\({\cal P}_{vz}\)
be the unique path in the tree~\({\cal T}_{n+1}\) with endvertices~\(v\) and~\(z.\)
If the path~\({\cal P}_{vz}\) does not contain~\(u\) as shown in Figure~\ref{mst_fig2}\((a),\)
then the edge~\((u,z)\) cannot be present in~\({\cal T}_{n+1}\) as this would create a cycle.
Removing the edge~\((u,v)\) and adding the edge~\((u,z),\) we get a new tree~\({\cal T}_{new}.\)
By construction, the sum of length of edges in~\({\cal T}_{new}\)
is strictly less than the sum of length of edges in the MST~\({\cal T}_{n+1},\)
a contradiction.

If the path~\({\cal P}_{vz}\) contains the node~\(u,\) then the edge~\((u,v)\)
necessarily belongs to~\({\cal P}_{vz}\) because~\((u,v)\) is the unique path in
the tree~\({\cal T}_{n+1}\) connecting~\(u\) and~\(v.\) In this case,
the edge~\((v,z)\) cannot be in~\({\cal T}_{n+1}\)
as this would create a cycle (see Figure~\ref{mst_fig2}\((b)\)). Define~\({\cal T}_{new}\) to be the graph
obtained by deleting the edge~\((u,v)\) and adding the edge~\((v,z).\)
The graph~\({\cal T}_{new}\) is again a tree and
the sum of length of edges in~\({\cal T}_{new}\)
is strictly less than the sum of length of edges in the MST~\({\cal T}_{n+1},\)
a contradiction.~\(\qed\)

%The above result allows us to estimate the difference in the length between MST's
%differing by one vertex. Recall from~(\ref{ave_per_sq}) that
%each~\(r_n \times r_n\) in square~\(\{S_l\}\) contains on an average between~\(2D_1 M\log{n}\)
%and~\(\frac{D_2}{2} M\log{n}\) vertices of~\(\{X_k\}_{1 \leq k \leq n}\) for some constants~\(D_1\) and~\(D_2.\)
%Here~\(M > 0\) is as in the definition of~\(r_n\) in~(\ref{rn_def_mst}).
%Adding one extra vertex~\(X_{n+1}\) we then have the following statement:
%On an average, each~\(r_n \times r_n\) square in~\(\{S_l\}\) contains between~\(2D_1M \log{n}\)
%and~\(\frac{D_2}{4} M \log{n} +1 \leq \frac{D_2}{2} M\log{n}\) vertices of~\(\{X_k\}_{1 \leq k \leq n+1}\) for all~\(n\) large.

\emph{Proof of Lemma~\ref{one_point_est}}: Suppose that the event~\(Z_{tot}(n+1)\) defined
in~(\ref{z_tot_def}) occurs and suppose the node~\(X_{j}  = v \in S_i\)
for some~\(1 \leq i \leq N.\)

To find an upper bound for~\(MST_{n+1} - MST_n(j),\)
let~\({\cal T}_n(j)\) be the MST containing
the nodes~\(\{X_k\}_{1 \leq k \neq j \leq n+1}.\) Since the event~\(Z_{j}(i) \supseteq Z_{tot}(n+1)\) occurs,
the square~\(S_i\) contains some node~\(w \in \{X_k\}_{ 1 \leq k \neq j \leq n}.\) Joining~\(v\)
and~\(w\) by an edge, we get a new tree containing all the nodes~\(\{X_k\}_{1 \leq k \leq n+1}.\)
The edge length between~\(v\) and~\(w\) is at most~\(r_n \sqrt{2}\) and so
\begin{equation}\label{mst_up_y}
MST_{n+1} - MST_n(j) \leq r_n \sqrt{2} \ind(Z_{tot}(n+1)).
\end{equation}

To obtain a lower bound for~\(MST_{n+1} - MST_n(j),\) we use property~\((h1)\)
and estimate the difference in length
of the MST obtained by removing the node~\(X_j = v\) from the MST~\({\cal T}_{n+1}\)
containing all the nodes~\(\{X_k\}_{1 \leq k \leq n+1}.\)
From property~\((h1),\) every edge in the MST~\({\cal T}_{n+1}\)
containing~\(v\) as an endvertex, has its other endvertex in some square~\(S_k \in {\cal N}_{20}(S_i).\)
Since~\(Z_{tot}(n+1) \) occurs, there are at most~\(4\epsilon_2 M \log{n}\) nodes
of~\(\{X_k\}_{1 \leq k \neq j \leq n}\) in every square~\(S_k \in {\cal N}_{20}(S_i)\)
(see definition of~\(Z_{tot}(n+1)\) prior to~(\ref{z_tot_def})).
There are at most~\(40^2\) squares of~\(\{S_k\}\) in~\({\cal N}_{20}(S_i)\) and so
the degree~\(d(v)\) of~\(v\) in the tree~\({\cal T}_{n+1}\) is at most
\begin{equation}\label{du_est}
d(v) \leq 40^2 \cdot (4\epsilon_2 M \log{n}) =  C_2 \log{n}
\end{equation}
for some constant~\(C_2  >0.\)

\begin{figure}
\centering
\begin{subfigure}{0.5\textwidth}
\centering
\includegraphics[width=2in, trim= 160 300 100 50, clip=true]{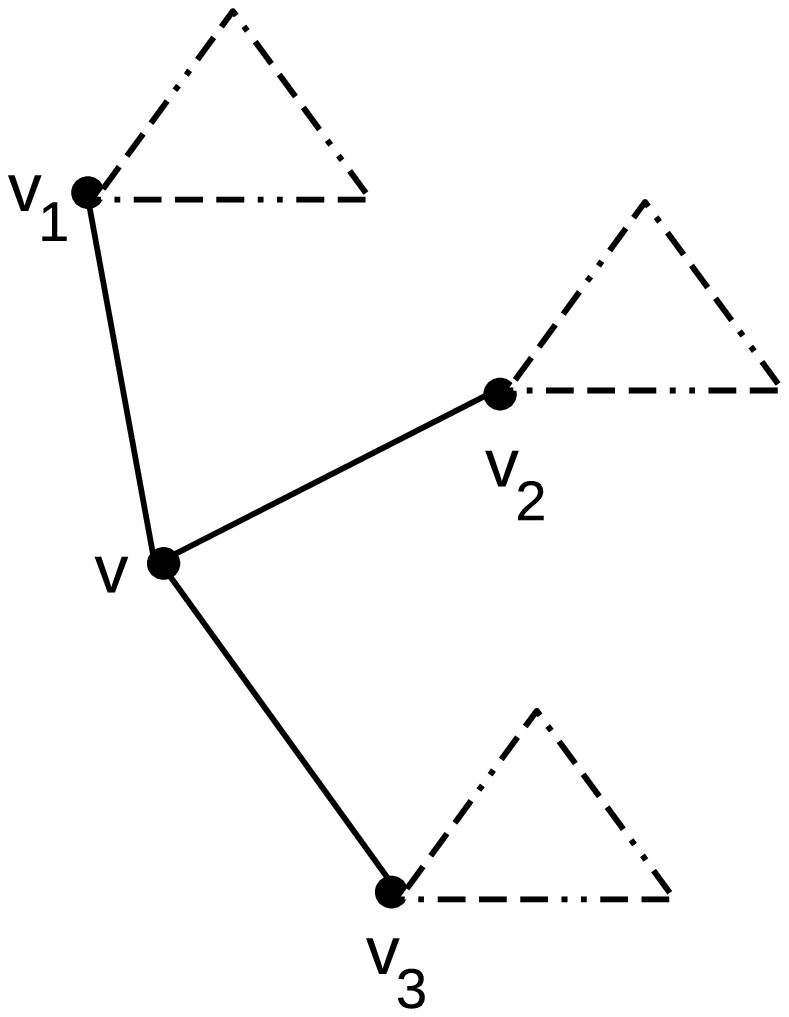}
   \caption{Before removing the node~\(X_j = v.\)}%\nonumber%{t1}
\end{subfigure}% seems this is important....
\begin{subfigure}{0.5\textwidth}
\centering
  \includegraphics[width=2in, trim= 160 300 100 50, clip=true]{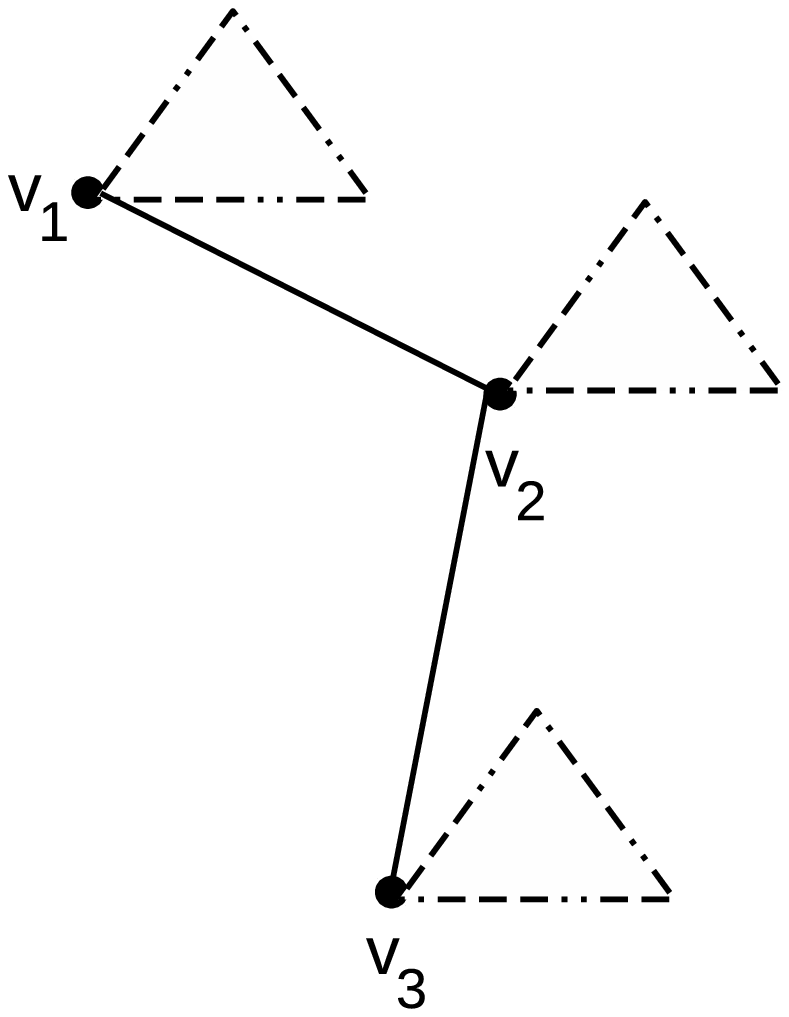}
   \caption{After removing~\(X_j = v.\)} %\nonumber%{t2}
  \end{subfigure}

\caption{Removing the vertex~\(X_{j} = v\) and forming a new tree.}
\label{mst_fig3}
\end{figure}

Suppose~\(\{v_{k}\}_{1 \leq k \leq d(v)}\) are the neighbours of~\(X_{j} = v\)
in the tree~\({\cal T}_{n+1}.\) Remove the node~\(v\) and the edges containing~\(v\) as an endvertex
and add the edges~\((v_k,v_{k+1})\) for~\(1 \leq k \leq d(v)-1\) as shown in Figure~\ref{mst_fig3}.
Here~\(d(v)= 3\) and the broken triangles represent the corresponding
subtrees of~\({\cal T}_{n+1}\) attached to the nodes~\(v_1,v_2\) and~\(v_3.\)

The resulting graph is a tree containing all the nodes~\(\{X_k\}_{1 \leq k \neq j \leq n+1}.\)
Each edge removed in the above process belongs to~\({\cal T}(n+1)\) and so has length
at most~\(20 r_n \sqrt{2}\) (property~\((h1)\)).
Using~(\ref{du_est}), the total length of the edges removed is then at most~
\[d(v) \cdot (20 r_n \sqrt{2}) \leq C_3 r_n \log{n}\] for some constant~\(C_3 > 0.\) Consequently
\begin{equation}\label{mst_down_y}
MST_n(j) \leq MST_{n+1} + C_3 r_n \log{n} \ind(Z_{tot}(n+1)).
\end{equation}
From~(\ref{mst_up_y}) and~(\ref{mst_down_y}), we obtain~(\ref{mst_diff})
for the case when~\(Z_{tot}(n+1)\) occurs.

If~\(Z_{tot}(n+1)\) does not occur, we use the crude upper bound
that any edge belonging to either of the spanning trees~\({\cal T}_{n+1}\) or~\({\cal T}_{n}(j)\)
has length most~\(\sqrt{2}\) and there are~\(n\) edges in~\({\cal T}_{n+1}\)
and~\(n-1 \leq n\) edges in~\({\cal T}_n.\) This proves~(\ref{mst_diff}).

To prove~(\ref{mst_diff_one}),
let~\(M > 0\) be large so that~\(\mathbb{P}(Z_{tot}(n+1)) \geq 1-\frac{1}{n^3}.\)
Setting~\(j = n\) in~(\ref{mst_diff}) and using the estimate for~\(r_n\) in~(\ref{rn_def_mst}), we then get
\begin{equation}\nonumber
\mathbb{E}|MST_{n+1} - MST_n| \leq C_2 \frac{(\log{n})^{3/2}}{\sqrt{n}} + n\sqrt{2} \frac{1}{n^3} \leq C_3\frac{(\log{n})^{3/2}}{\sqrt{n}}
\end{equation}
This proves~(\ref{mst_diff_one}).~\(\qed\)

%WRT MOTRE +etc...
%variance of the
%length of the minimum length cycle.\\

\emph{Proof of~\ref{var_mst_est_main} of Theorem~\ref{var_mst_thm}}:
We use the martingale difference method and
for~\(1 \leq j \leq n+1,\) let~\[{\cal F}_j = \sigma\left(X_1,\ldots,X_j\right)\] denote the sigma field
generated by the random variables~\(X_1,\ldots,X_i.\) Defining the martingale difference
\begin{equation}\label{gi_diff}
G_j = \mathbb{E}(MST_{n+1} | {\cal F}_j) - \mathbb{E}(MST_{n+1} | {\cal F}_{j-1}),
\end{equation}
we have that~\[MST_{n+1} -\mathbb{E}MST_{n+1} = \sum_{j=1}^{n+1} G_j\] and so by the martingale property
\begin{equation}\label{var_tsp_gi}
var(MST_{n+1})  = \left(\sum_{j=1}^{n+1} G_j\right)^2 = \sum_{j=1}^{n+1} \mathbb{E}G_j^2.
\end{equation}
There is a constant~\(C> 0\) such that
\begin{equation}\label{gi_sec}
\max_{1 \leq j \leq n+1}\mathbb{E}G_j^2 \leq \frac{C(\log{n})^3}{n}
\end{equation}
for all~\(n \geq 1\) and this proves~(\ref{var_mst_est_main}).

To prove~(\ref{gi_sec}), we rewrite~\(G_j\) in a more convenient form.
Let~\(\omega = (x_1,\ldots,x_{n+1})\) and~\(\omega' = (y_1,\ldots,y_{n+1})\)
be two vectors in~\((\mathbb{R}^2)^{n+1}.\) We say that~\(\{x_k\}_{1 \leq k \leq n+1}\)
are the nodes of~\(\omega.\)  Defining~\(\omega_j = (x_1,\ldots,x_j,y_{j+1},\ldots,y_{n+1})\) for~\(1 \leq j \leq n+1\)
and using Fubini's theorem, we get
\begin{equation}
|G_j| = \left|\int (M(\omega_j) - M(\omega_{j-1}))f(y_i)\ldots f(y_{n+1})dy_{j}\ldots dy_{n+1} \right|  \leq H_j,
\end{equation}
where
\begin{equation}
H_j := \int |M(\omega_j) - M(\omega_{j-1})|f(y_j)\ldots f(y_{n+1})dy_{j}\ldots dy_{n+1}, \label{hj_def}
\end{equation}
and~\(M(\omega_j)\) is the length of the MST containing all the nodes in~\(\omega_i.\)\\
\emph{Proof of~(\ref{gi_sec})}: Let~\(Z_{tot}(n+1)\) be the event defined in~(\ref{z_tot_def})
prior to the proof of property~\((h2)\) above.
From~(\ref{hj_def}),
\begin{equation}\label{g_decom}
H_j = I_1  + I_2,
\end{equation}
where
\begin{eqnarray}
I_1 &=& \int |M(\omega_j) - M(\omega_{j-1})|\ind(\omega_j \in Z_{tot}(n+1))\ind(\omega_{j-1} \in Z_{tot}(n+1)) \nonumber\\
&&\;\;\;\;\;\;\;f(y_j)\ldots f(y_{n+1})dy_{j}\ldots dy_{n+1} \label{i1_def_f2}
\end{eqnarray}
and~\(I_2 = I_1 - H_j.\)

We have that
\begin{equation}\label{i_ests}
\mathbb{E}I_1^2 \leq \frac{C(\log{n})^3}{n} \text{ and } \mathbb{E}I_2^2 \leq \frac{4}{n}
\end{equation}
for some constant~\(C > 0\) and all~\(n\) large. Since~\(|G_j|^2 \leq H_j^2 = (I_1+I_2)^2 \leq 2(I_1^2 + I_2^2),\)
we get that \[\mathbb{E}(G_j^2) \leq 2 \left(\frac{C(\log{n})^3}{n} + \frac{4}{n}\right) \leq \frac{3C(\log{n})^3}{n},\]
proving~(\ref{gi_sec}).

We obtain the estimates for~\(\mathbb{E}I^2_1\) and~\(\mathbb{E}I^2_2\) in~(\ref{i_ests}), separately below.\\
\underline{Estimate for~\(I_1\)}: Let~\({\cal T}_n(j)\) be the MST containing
all the vertices~\(\{x_k\}_{1 \leq k \leq j-1} \cup \{y_k\}_{j+1 \leq k \leq n}.\)
If~\(L({\cal T}_n(j))\) is the length of~\({\cal T}_n(j),\) then from~(\ref{mst_diff})
we have for~\(t \in \{j-1,j\}\) that
\begin{equation}
|M(\omega_t) - L({\cal T}_n(j))| \ind(\omega_t \in Z_{tot}(n+1)) \leq C r_n\log{n} \label{t_om_td2}
\end{equation}
for some constant~\(C > 0.\) From~(\ref{t_om_td2}),~(\ref{i1_def_f2}) and triangle inequality, we therefore have
\begin{equation}\label{i1_estr2}
I_1 \leq  2C r_n (\log{n}) \text{   and so   } \mathbb{E}(I^2_1) \leq 4C^2 r_n^2 (\log{n})^2  \leq C_1 \frac{(\log{n})^3}{n}
\end{equation}
for some constant~\(C_1  >0.\) The final estimate in~(\ref{i1_estr2}) follows from the expression for~\(r_n\)
in~(\ref{rn_def_mst}).

\underline{Estimate for~\(I_2\)}: To estimate~\(I_2,\) use the fact the MST containing
all the nodes of~\(\omega_t, t = j-1,j\) has~\(n\) edges, each of which has length at most~\(\sqrt{2}.\) Therefore
\begin{eqnarray}
I_2 &\leq& \int n\sqrt{2} \left(\ind(\omega_j \notin Z_{tot}(n+1)) + \ind(\omega_{j-1} \notin Z_{tot}(n+1) \right) \nonumber\\
&&\;\;\;\;\;\;\;f(y_j)\ldots f(y_{n+1})dy_j\ldots dy_{n+1} \nonumber\\
&& = J_1 + J_2,
\end{eqnarray}
where~\(J_1 = n \sqrt{2} \int \ind(\omega_j \notin Z_{tot}(n+1))f(y_j)\ldots f(y_n)dy_j\ldots dy_n\) and~\(J_2\) is the remaining term.
Using Cauchy-Schwarz inequality,
\[J_1^2 \leq 2n^2 \left(\mathbb{E}(\ind(Z^c_{tot}(n+1))|{\cal F}_j)\right)^2 \leq 2n^2 \mathbb{E}(\ind(Z^c_{tot}(n+1))|{\cal F}_j)\]
Similarly~\(J_2^2 \leq 2n^2 \mathbb{E}\left(\ind(Z^c_{tot}(n+1)) | {\cal F}_{j-1}\right).\)
Using~\(I_2^2 \leq 2(J_1^2 + J_2^2)\)
and the fact that~\(\mathbb{E}(\mathbb{E}(X|{\cal F}_j)|{\cal F}_{j-1}) = \mathbb{E}(X | {\cal F}_{j-1}),\)
we get~\[\mathbb{E}(J_1^2 + J_2^2|{\cal F}_{j-1}) \leq 4n^2 \mathbb{P}(Z_{tot}^c(n+1) | {\cal F}_{j-1}).\]

Since~\(I_2^2 \leq (J_1+J_2)^2 \leq 2(J_1^2 + J_2^2),\)
we get \[\mathbb{E}(I_2^2) \leq 4n^2 \mathbb{P}(Z_{tot}(n+1)^c) \leq 4n^2 e^{-CM\log{n}}\] for some constant~\(C>  0,\) using~(\ref{z_tot_est}).
Letting~\(M >0\) large so that~\(e^{-CM\log{n}} \leq \frac{1}{n^3},\) we get the estimate for~\(I_2\) in~(\ref{i_ests}).~\(\qed\)

Using the variance estimate~(\ref{var_mst_est_main}), we prove the almost sure convergence result.\\
\emph{Proof of~(\ref{as_conv_mst}) in Theorem~\ref{var_mst_thm}}:
From~(\ref{var_mst_est_main}) and Borel-Cantelli lemma,
\begin{equation}\label{sub_conv}
\frac{1}{n}\left(MST_{n^2} - \mathbb{E} MST_{n^2}\right) \longrightarrow 0 \text{ a.s. }
\end{equation}
as~\(n \rightarrow \infty.\)
For convergence along the sequence~\(a_n = n,\) we use a subsequence argument and define
\begin{equation}\label{d_def_mst}
D_{n} := \max_{n^2 \leq k < (n+1)^2} \left|MST_{k} - MST_{n^2}\right|.
\end{equation}
Recalling the event~\(Z_{tot}(n+1)\) defined in~(\ref{z_tot_def}), let
\begin{equation}\label{y_tot_def}
Y_{tot}(n) := \bigcap_{n^2 \leq k < (n+1)^2} Z_{tot}(k+1)
\end{equation}
so that from~(\ref{mst_diff}), the difference
\begin{equation}
|MST_{k+1} - MST_{k}| \leq C_1 r_{k} (\log{k}) \ind(Y_{tot}(n)) + k \sqrt{2} \ind(Y_{tot}^c(n)) \nonumber
\label{k_dif}
\end{equation}
for each~\(n^2 \leq k <(n+1)^2\) and for some constants~\(C_1,C_2 > 0\) not depending on~\(k\) or~\(n.\)

%\begin{eqnarray}
%|MST_{k+1} - MST_{k}| \leq C_1 r_{k} (\log{k}) \ind(Y_{tot}(n)) + k \sqrt{2} \ind(Y_{tot}^c(n)) \nonumber\\
%&\leq& C_2 \frac{(\log{k})^{3/2}}{\sqrt{k}} \ind(Y_{tot}(n)) + (n+1)^2 \sqrt{2} \ind(Y_{tot}^c(n)) \nonumber\\
%\label{k_dif}
%\end{eqnarray}
%for each~\(n^2 \leq k <(n+1)^2\) and for some constants~\(C_1,C_2 > 0\) not depending on~\(k\) or~\(n.\)
%The final estimate~(\ref{k_dif}) is true since from~(\ref{rn_def_mst}) we get that~\(r_k \leq C_3 \sqrt{\frac{\log{k}}{k}}\)
%for some constant~\(C_3 > 0.\)

From~(\ref{rn_def_mst}))
we have that~\(r_k \leq C_2 \sqrt{\frac{\log{k}}{k}} \leq C_3 \frac{\sqrt{\log{n}}}{n}\)
for some positive constants~\(C_2,C_3\) and so
\begin{equation}\label{rk_est}
r_k \log{k} \leq C_3 \frac{(\log{n})^{3/2}}{n} \text{ and }k\sqrt{2} \leq (n+1)^2 \sqrt{2}
\end{equation}
for some positive constants~\(C_2,C_3,C_4\) and for all~\(n^2 \leq k < (n+1)^2.\)
Using~(\ref{rk_est}) in~(\ref{k_dif}) and adding telescopically, we get
\begin{equation}\label{sub_dif}
|MST_{k} - MST_{n^2}| \leq C_4 \frac{(\log{n})^{3/2}}{n} (k-n^2) \ind(Y_{tot}(n)) + (k-n^2)(n+1)^2 \sqrt{2} \ind(Y_{tot}^c(n))
\end{equation}
for~\(n^2 \leq k < (n+1)^2.\)

From~(\ref{d_def_mst}),~(\ref{sub_dif}) and the fact that~\(k-n^2 \leq (n+1)^2 - n^2 \leq 4n\) for all~\(n\) large, we get
\begin{equation}\label{d_est_ac}
D_{n} \leq C_5 \left(\log{n}\right)^{3/2} \ind(Y_{tot}(n)) + 4n(n+1)^2\ind(Y_{tot}^c(n)).
\end{equation}
From the estimate for~\(Z_{tot}(k)\) in~(\ref{z_tot_est})
\begin{eqnarray}
\mathbb{P}(Y_{tot}(n)) &\geq& 1 - \sum_{k=n^2}^{(n+1)^2-1} \mathbb{P}(Z_{tot}^c(k)) \nonumber\\
&\geq& 1-\sum_{k=n^2}^{(n+1)^2-1} \exp\left(-C M\log{k}\right) \nonumber\\
&\geq& 1-((n+1)^2-n^2)\exp\left(-C M\log(n^2)\right), \nonumber
\end{eqnarray}
for all~\(n\) large and for some constant~\(C >0.\) Setting~\(M > 0\) large so that~\(\exp\left(-CM\log(n^2)\right) \leq \frac{1}{n^{10}}\) we then get that
\begin{equation}
\mathbb{P}(Y_{tot}(n)) \geq 1-\frac{(2n+1)}{n^{9}} \geq 1-\frac{1}{n^7} \label{y_tot_est}
\end{equation}
for all~\(n\) large.

From Borel-Cantelli lemma and~(\ref{y_tot_est}) we get that~\(\mathbb{P}(\liminf_n Y_{tot}(n) ) = 1\) and so a.s.\ \(\ind(Y_{tot}^c(n)) = 0\) for all large~\(n.\) From~(\ref{d_est_ac}) and~(\ref{y_tot_est}), we therefore get
\begin{equation}\label{dn_as}
\frac{D_n}{n} \leq \frac{C_5(\log{n})^{3/2}}{n} + 4n(n+1)^2\ind(Y_{tot}^c(n)) \longrightarrow 0 \text{ a.s.\ }
\end{equation}
and
\begin{equation}\label{dn_exp}
\frac{\mathbb{E}D_n}{n} \leq \frac{C_5(\log{n})^{3/2}}{n} + \frac{4n(n+1)^2}{n^7} \longrightarrow 0
\end{equation}
as~\(n \rightarrow \infty.\)

Finally for~\(n^2 \leq k <(n+1)^2,\) write
\begin{eqnarray}
\frac{1}{\sqrt{k}}\left|MST_k - \mathbb{E}MST_k\right| &\leq& \frac{1}{\sqrt{k}} |MST_k - MST_{n^2}|
+ \frac{1}{\sqrt{k}} \mathbb{E}|MST_k - MST_{n^2}| \nonumber\\
&\leq& \frac{1}{n} |MST_k - MST_{n^2}| + \frac{1}{n} \mathbb{E}|MST_k - MST_{n^2}| \nonumber\\
&\leq& \frac{D_n}{n} + \frac{\mathbb{E}D_n}{n} \nonumber
\end{eqnarray}
and use~(\ref{dn_as}) and~(\ref{dn_exp}) to get that~\(\frac{1}{\sqrt{k}}\left(MST_k - \mathbb{E}MST_k\right)\)
converges to zero a.s.\ as~\(k \rightarrow \infty.\) ~\(\qed\)

%SEE ABV +eTC...
%ALSO SHOW IN EXP...

%Arguing iteratively as in the proof of~\((h2),\) we get
%\begin{eqnarray}
%\left|MST_{k} - MST_{n^2}\right| &\leq& C_1r_{n^2} (\log{n^2})(k-n^2)\ind\left(Y_{tot}(n^2)\right) \nonumber\\
%&&\;\;\;\;\;\;+\;\;\;(k-n^2) n^2\sqrt{2} \ind\left(Y_{tot}^c(n^2)\right)\label{d_est_mst}
%\end{eqnarray}
%GOOTHALS HERE...

\emph{Proof of~(\ref{exp_mst_u}) and~(\ref{eq_mst1_u}) in Theorem~\ref{var_mst_thm}}: The variance estimate~(\ref{var_mst_est_main}) is proved above.
The upper bound for~\(\mathbb{E} MST_n\) in~(\ref{exp_mst_u})
is obtained from the strips estimate~(\ref{mst_ab}) with~\(a = n\) and~\(b = 1.\) This also proves~(\ref{eq_mst1_u}).

To prove the lower bound for~\(\mathbb{E}MST_n\) in~(\ref{exp_mst_u}), let~\(l(X_{j},{\cal T}_n)\) denote the total
length of the edges  containing the node~\(X_j\) in the MST~\({\cal T}_n.\) From~(\ref{len_cyc_def}),
\(MST_n = \frac{1}{2} \sum_{j=1}^{n} l(X_j,{\cal T}_n) \geq \frac{1}{2} \sum_{j=1}^{n} d(X_j,\{X_k\}_{k \neq j}),\)
where~\(d(X_j,\{X_k\}_{k \neq j})\) is the minimum distance of the node~\(X_j\) from all the other nodes.
Therefore\\\(\mathbb{E} MST_n \geq \frac{n}{2} \mathbb{E}d(X_1,\{X_j\}_{2 \leq j \leq n}) \geq C_1\sqrt{n}\) for some constant~\(C_1 > 0,\)
by arguing analogous to the proof of~(\ref{x2_est}) in property~\((b2).\)~\(\qed\)

%at the beginning of Section~\ref{pf_mst}).
%Dividing the unit square into strips of size~\(\frac{1}{\sqrt{n}} \times 1\) and arguing as in the proof of~\((b4),\)
%we obtain~\(MST_n \leq \sqrt{n} + n\frac{1}{\sqrt{n}} \sqrt{2} + 2 \leq 3\sqrt{n}\) for all~\(n\) large.

%\otimes_{k=1}^{\infty} \Omega_k

\emph{Proof of~(\ref{eq_mst2_u}) in Theorem~\ref{var_mst_thm}}:
We perform Poissonization and construct a Poisson process~\({\cal P}\)
in the unit square~\(S\) with intensity~\(nf(.)\) as follows.
Let~\(\{V_{i,k}\}_{1 \leq i \leq N, k\geq 1}\) be i.i.d.\ random vectors in~\(\mathbb{R}^2\)
with density\\\(\frac{f(x)}{\int_{S_i}f(x)dx}\ind(x \in S_i).\)
Let~\(\{N(S_i)\}_{1 \leq i \leq N}\) be independent Poisson random variables
such that~\(N(S_i)\) has mean~\(n\int_{S_i} f(x)dx\) for~\(1 \leq i \leq T.\)
The random variables ~\(\{N(S_i)\}\) are independent of~\(\{V_{i,k}\}\)
and we define~\((\{V_{i,k}\}, \{N(S_i)\})\) on the probability space~\((\Omega_0, {\cal F}_0, \mathbb{P}_0).\)

For~\(1 \leq i \leq N,\) if~\(N(S_i) \geq 1,\) then we set~\(\{V_{i,k}\}_{1 \leq k \leq N(S_i)}\) to be the nodes of~\({\cal P}\) in the square~\(S_i.\)
Analogous to~(\ref{min_weight_tree}), let~\({\cal T}^{(P)}_n\) be
the MST containing all the nodes of~\({\cal P}\)
in the unit square~\(S\) and as in~(\ref{min_weight_tree}) define~\(MST^{(P)}_n := L({\cal T}^{(P)}_n).\)

We find lower bounds for the length~\(MST^{(P)}_n\) in the Poisson process and then later convert
the estimates to the Binomial process.
We first need some preliminary definitions and computations. Analogous to~(\ref{ave_per_sq}), we have for every~\(1 \leq i \leq N\) that
\begin{equation} \label{ave_per_sq_poi}
2\epsilon_1 M \log{n} \leq n\epsilon_1 r_n^2 \leq \mathbb{E}_0 N(S_i) \leq n\epsilon_2 r_n^2  \leq 3\epsilon_2 M \log{n}
\end{equation}
where~\(\epsilon_1,\epsilon_2 > 0\) is as in~(\ref{f_eq}). Defining
\begin{equation}\label{yj_def}
Y_i := \left\{\epsilon_1 M\log{n} \leq N(S_i) \leq 4\epsilon_2 M\log{n}\right\}
\end{equation}
we get by standard Poisson distribution estimates (Theorem~A.1.15, pp. 313, Alon and Spencer (2008)) that
\begin{equation}\label{yk_est}
\mathbb{P}_0(Y_i) \geq 1-e^{-2CM\log{n}}
\end{equation}
for some constant~\(C > 0\) not depending on~\(M\) and for all~\(n\) large.

For~\(q \geq 1,\) recall the definition of the~\(q-\)neighbourhood~\({\cal N}_{q}(S_i)\)
of the square~\(S_i, 1 \leq i \leq N,\) from the discussion following Lemma~\ref{one_point_est}.
Let~\(W_1,\ldots,W_{T}\) be a maximal set of squares in~\(\{S_k\}\) such that~\({\cal N}_{20}(W_i) \bigcap {\cal N}_{20}(W_j) = \emptyset\)
for any~\(1 \leq i\neq j \leq T.\) There are~\((41)^2\) squares in~\({\cal N}_{20}(W_i)\) for any square~\(W_i\) and so by our choice of~\(r_n\) in~(\ref{rn_def_mst}), we have that~\(\bigcup_{i=1}^{T} {\cal N}_{20}(W_i) = \bigcup_{k=1}^{N} S_k.\)
Since there are a total~\(N = \left(\frac{1}{r_n}\right)^2\) squares in~\(\{S_k\},\) we must have
\begin{equation}\label{t_est}
C_1\frac{n}{\log{n}} \leq T = \frac{N}{(41)^2} = \left(\frac{1}{41r_n}\right)^2 \leq C_2\frac{n}{\log{n}}
\end{equation}
for some positive constants~\(C_1,C_2,\) using the bounds for~\(r_n\) in~(\ref{rn_def_mst}).
For~\(1 \leq i \leq T,\) let
\begin{equation}\label{qi_def}
Q_i := \bigcap_{k : S_k \in {\cal N}_{20}(W_i)} Y_k,
\end{equation}
so that from~(\ref{yk_est}) we get
\begin{equation}\label{qi_est}
\mathbb{P}_0(Q_i) \geq 1-(41)^2e^{-2CM\log{n}}
\end{equation}
for some constant~\(C > 0.\)

The event~\(Q_i\) is useful in the following way.\\
\((h2)\) Suppose the event~\(Q_i\) occurs for some~\(1 \leq i \leq T\)
and let~\(e\) be an edge of the MST~\({\cal T}^{(P)}_n\) containing a node~\(v \in W_i.\)
If~\(u\) denotes the other endvertex of~\(e,\) then~\(u \in S_k\) for some~\(S_k \in {\cal N}_{20}(W_i).\)\\\\
The proof of~\((h2)\) is analogous to the proof of property~\((h1)\) stated below Lemma~\ref{one_point_est}.\\

Recall from paragraph prior to~(\ref{ave_per_sq_poi}) that~\(N(W_i)\) is the number of nodes of the Poisson process~\({\cal P}\)
in the square~\(W_i\) and that~\(\{V_{i,k}\}_{1 \leq k \leq N(W_i)}\) are  the nodes of~\({\cal P}\) in~\(W_i.\) Let~\(l(V_{i,k}, {\cal T}^{(P)}_n)\) be the sum of length of the edges containing the node~\(V_{i,k}\) as an endvertex
in the MST~\({\cal T}^{(P)}_n,\) with the notation that the sum length is zero if~\(N(W_i) = 0.\) From~(\ref{len_cyc_def})
\(MST^{(P)}_n = L({\cal T}^{(P)}_n)\) satisfies
\begin{equation}\label{mst_p_est}
MST^{(P)}_n \geq \frac{1}{2} \sum_{i=1}^{T} \sum_{k=1}^{N(W_i)} l\left(V_{i,k}, {\cal T}^{(P)}_n\right) \geq \frac{1}{2} \sum_{i=1}^{T} \sum_{k=1}^{N(W_i)} l\left(V_{i,k}, {\cal T}^{(P)}_n\right)\ind(Q_i).
\end{equation}
If the event~\(Q_i\) occurs, the number of nodes~\(N(W_i) \geq \epsilon_1 M\log{n}.\) Moreover, from property~\((h2)\) above,
every edge containing~\(V_{i,k} \in W_i\) as an endvertex has its other endvertex in some square belonging to the neighbourhood~\({\cal N}_{20}(W_i).\)
Therefore~\[l\left(V_{i,k}, {\cal T}^{(P)}_n\right)\ind(Q_i) \geq d_{i,k} \ind(Q_i)\]
where~\(d_{i,k}\) is the minimum distance of the node~\(V_{i,k} \in W_i\) from all the nodes of~\({\cal P}\) in~\({\cal N}_{20}(W_i).\)

Summarizing,
\begin{equation}\label{mst_p_est2}
MST^{(P)}_n  \geq \sum_{k=1}^{\epsilon_1 M\log{n}} \sum_{i=1}^{T} F_{i,k},
\end{equation}
where~\(F_{i,k} := d_{i,k} \ind(Q_i).\) We need the following property regarding the moments of~\(F_{i,k}.\)\\
\((h3)\) There are positive constants~\(C_1,C_2\) and~\(C_3\) such that for any~\(1 \leq i \leq T\)
and any~\(1 \leq k \leq M\log{n},\)
\begin{equation}\label{fjk_est_mean}
C_1 \frac{r_n}{\sqrt{\log{n}}} \leq \mathbb{E}_0 F_{i,k} \leq C_2 \frac{r_n}{ \sqrt{\log{n}}}
\text{ and } \mathbb{E}_0 F^2_{i,k} \leq C_3 \frac{r_n^2}{\log{n}}.
\end{equation}
\emph{Proof of~\((h3)\)}: There are~\(L = (41)^2\)~
squares of~\(\{S_k\}\) in~\({\cal N}_{20}(W_i)\) and if the event~\(Q_i\) occurs,
then each square~\(S_k \in {\cal N}_{20}(W_i)\)
has between~\(\epsilon_1 M\log{n}\) and~\(4\epsilon_2 M\log{n}\) nodes of~\({\cal P}\) (see~(\ref{qi_def}) and~(\ref{yj_def})).

For positive integers~\(l_1,\ldots,l_L\) define
\[E(l_1,\ldots,l_{L}) = \bigcap_{S_k \in {\cal N}_{20}(W_i)} \{N(S_k) = l_k\}\]
and use the definition of~\(Q_i\) in~(\ref{qi_def})
to get that~\(Q_i = \bigcup_{(l_1,\ldots,l_L)} E(l_1,\ldots,l_L),\)
where the union is over all~\(L-\)tuples satisfying
\begin{equation}\label{l_sat}
\epsilon_1 M\log{n} \leq l_k \leq 4\epsilon_2 M \log{n}, 1 \leq k \leq L.
\end{equation}
If~(\ref{l_sat}) holds, then arguing as in the proof of~(\ref{x2_est}) in property~\((b2),\) we get
\begin{equation}\nonumber
C_4 \frac{r_n}{\sqrt{\log{n}}} \leq \mathbb{E}_0 \left(d_{i,k}\left|\right.E(l_1,\ldots,l_{L})\right) \leq C_5 \frac{r_n}{\sqrt{\log{n}}}
\end{equation}
and~\(\mathbb{E}_0 \left(d^2_{i,k}\left|\right.E(l_1,\ldots,l_{L})\right) \leq C_6 \frac{r^2_n}{{\log{n}}}\)
for some positive constants~\(C_4,C_5\) and~\(C_6,\) not depending on~\(\{l_k\}\) or~\(i.\)
Thus
\begin{equation} \nonumber
C_4 \frac{r_n}{\sqrt{\log{n}}} \mathbb{P}_0(Q_i) \leq
\mathbb{E}_0 \left(d_{i,k}\ind(Q_i)\right) \leq C_5 \frac{r_n}{\sqrt{\log{n}}} \mathbb{P}_0(Q_i).
\end{equation}
and~\(\mathbb{E}_0 \left(d^2_{i,k}\ind(Q_i)\right) \leq C_6 \frac{r^2_n}{\log{n}} \mathbb{P}_0(Q_i).\)
Using the estimate for~\(\mathbb{P}_0(Q_i)\) in~(\ref{qi_est}) we then get~(\ref{fjk_est_mean}).

%The total number of points~\(N_{tot}(j)\) in all the squares of~\({\cal N}_{20}(S_j)\)
%satisfies~\(C_4\log{n} \leq N_{tot}(j) \leq C_5 \log{n}\) for some positive constants~\(C_4,C_5.\)

From~(\ref{fjk_est_mean}) and the Paley-Zygmund inequality~(\ref{paley}), we have for~\(1 \leq i \leq T\) and~\(1 \leq k \leq \epsilon_1 M\log{n}\) that
\begin{equation}\label{play}
\mathbb{P}_0\left( F_{i,k} \geq \delta_1 \frac{r_n}{\sqrt{\log{n}}}\right) \geq \delta_2
\end{equation}
for some positive constants~\(\delta_1\) and~\(\delta_2,\) not depending on~\(i\) or~\(k.\)
We use~(\ref{play}) to lower bound~\(MST^{(P)}_n\)
in~(\ref{mst_p_est2}) as follows. Let~\(G_{i,k} = \{F_{i,k} \geq \delta_1 \frac{r_n}{\sqrt{\log{n}}}\}\) and use~(\ref{mst_p_est2})
to get
\begin{eqnarray}\label{mst_p_est3}
MST^{(P)}_n  \geq \sum_{k=1}^{\epsilon_1 M\log{n}} \sum_{i=1}^{T} F_{i,k} \ind(G_{i,k}) \geq \delta_1 \frac{r_n}{\sqrt{\log{n}}}\sum_{k=1}^{\epsilon_1 M\log{n}} \sum_{i=1}^{T} \ind(G_{i,k}). \nonumber
\end{eqnarray}

Since the Poisson process is independent on disjoint sets, the terms~\(F_{i_1,k}\) and~\(F_{i_2,k}\) are independent
for distinct~\(1 \leq i_1 \neq i_2 \leq T.\) Therefore we get from~(\ref{play}) and standard Chernoff estimates for Bernoulli random variables
that
\begin{equation}\label{cher_est}
\mathbb{P}_0\left(\sum_{i=1}^{T} \ind(G_{i,k}) \geq \delta_3 T\right) \geq 1-e^{-\delta_4 T}
\end{equation}
for some positive constants~\(\delta_3,\delta_4\) not depending on~\(k.\) Using the bounds for~\(T\) in~(\ref{t_est}), we get
\begin{equation}\label{cher_est2}
\mathbb{P}_0\left(\sum_{i=1}^{T} \ind(G_{i,k}) \geq \delta_5 \frac{n}{\log{n}}\right) \geq 1-\exp\left(-\delta_6 \frac{n}{\log{n}}\right)
\end{equation}
for some positive constants~\(\delta_5,\delta_6.\) Consequently,
\begin{eqnarray}\label{g_tot_est}
&&\mathbb{P}_0\left(\sum_{k=1}^{\epsilon_1 M\log{n}}\sum_{i=1}^{T} \ind(G_{i,k})  \geq  \left(\delta_5 \frac{n}{\log{n}}\right) \cdot \epsilon_1 M\log{n}\right) \nonumber\\
&&\;\;\;\;\geq 1-(\epsilon_1 M\log{n})\exp\left(-\delta_6 \frac{n}{\log{n}}\right) \nonumber\\
&&\;\;\;\; \geq 1-\exp\left(-\delta_7 \frac{n}{\log{n}}\right)
\end{eqnarray}
for all~\(n\) large, for some constant~\(\delta_7 > 0.\)

Using~(\ref{g_tot_est}) in~(\ref{mst_p_est3}) we get that with~\(\mathbb{P}_0-\)probability at least\\\(1-\exp\left(-\delta_7 \frac{n}{\log{n}}\right),\) the term
\begin{equation}\label{mst_p_est4}
MST^{(P)}_n  \geq \delta_1 \frac{r_n}{\sqrt{\log{n}}} \left(\delta_5 \frac{n}{\log{n}}\right) \cdot \epsilon_1 M\log{n} \geq C \sqrt{n},
\end{equation}
for some constant~\(C > 0,\) using the lower bound~\(r_n \geq \sqrt{\frac{M\log{n}}{n}}\) from~(\ref{rn_def_mst}).

%SEE FROM BLW...

Finally, to convert the estimates to the length~\(MST_n\) of the MST in the Binomial process, we let
\begin{equation} \nonumber
A := \{MST_n \geq C \sqrt{n}\}, A_P = \{MST^{(P)}_n \geq C\sqrt{n}\}
\end{equation}
and use dePoissonization formula~\(\mathbb{P}(A) \geq 1-D\mathbb{P}_0(A^c_P)\sqrt{n}\) for some constant~\(D > 0\) (see~(\ref{de_poiss_ax})). From~(\ref{mst_p_est4}) we then get~(\ref{eq_mst2_u}).~\(\qed\)

\emph{Proof of~(\ref{beta_conv})}: We need some preliminary definitions and estimates. For a set of nodes~\(x_1,\ldots,x_n\) in the unit square~\(S,\)
recall from Section~\ref{intro} that\\\(K_n(x_1,\ldots,x_n)\)
is the complete graph formed by joining all the nodes by straight line segments and~\(MST(x_1,\ldots,x_n)\)
is the length of the minimum spanning tree of~\(K_n(x_1,\ldots,x_n).\)

For any~\(a > 0,\) consider the graph~\(K_n(ax_1,\ldots,ax_n)\) where the length of the edge between
the vertices~\(ax_1\) and~\(ax_2\) is simply~\(a\) times the length of the edge between~\(x_1\) and~\(x_2\)
in the graph~\(K_n(x_1,\ldots,x_n)\)
Using the definition of MST in~(\ref{min_weight_tree}) we then have
\begin{equation}\label{mst_sc1}
MST(ax_1,\ldots,ax_n) = aMST(x_1,\ldots,x_n).
\end{equation}
Therefore if~\(Y_1,\ldots,Y_n\) are~\(n\) nodes uniformly distributed in the square~\(aS\) of side length~\(a,\)
then we get from~(\ref{mst_sc1}) that
\[MST(n;a) := MST(Y_1,\ldots,Y_n) = aMST(X_1,\ldots,X_n),\]
where~\(X_i = \frac{Y_i}{a}, 1 \leq i \leq n\) are i.i.d.\ uniformly distributed in~\(S.\)
Recalling the notation~\(MST_n = MST(X_1,\ldots,X_n)\) (see paragraph prior to Theorem~\ref{var_mst_thm}) we therefore get
\begin{equation}\label{mst_scale}
\mathbb{E} MST(n;a) = a\mathbb{E}MST_n.
\end{equation}

The following property is also needed for future use.\\
\((t1)\) For any positive integers~\(n_1,n_2 \geq 1\) we have that
\begin{eqnarray}
MST_{n_1+n_2} \leq MST_{n_1} + 3\sqrt{n_2} + \sqrt{2}. \label{mst_ab2}
\end{eqnarray}
\emph{Proof of~\((t1)\)}:
Let~\({\cal T}_1\) be the MST formed by the~\(n_1\) nodes~\(\{X_i\}_{1 \leq i \leq n_1}\)
and let~\({\cal T}_2\) be the MST formed by the remaining~\(n_2\) nodes. Joining~\({\cal T}_1\)
and~\({\cal T}_2\) by an edge~\(e_{12},\) we get a tree containing all the~\(n_1+n_2\) nodes.
Since~\(e_{12}\) has length at most~\(\sqrt{2},\) we get
\begin{equation}\label{mst_ab3}
MST_{n_1+n_2} \leq MST_{n_1} + MST(X_{n_1+1},\ldots,X_{n_2}) + \sqrt{2}.
\end{equation}
Using the strips estimate~(\ref{mst_ab}), the middle term in~(\ref{mst_ab3})
is bounded above by~\(3n_2\sqrt{2}.\)~\(\qed\)

To prove~(\ref{beta_conv}), it suffices to see that
\begin{equation}\label{sub_mst}
\frac{\mathbb{E}MST_{n^2}}{n} \longrightarrow \beta
\end{equation}
as~\(n \rightarrow \infty\) for some constant~\(\beta > 0.\) To see this is true,
use the definition of~\(D_n = \max_{n^2 \leq k < (n+1)^2} |MST_{k} - MST_{n^2}|\) in~(\ref{d_def_mst}) to get for~\(n^2 \leq k < (n+1)^2\) that
\[\frac{\mathbb{E}MST_k}{\sqrt{k}} \leq \frac{\mathbb{E}MST_k}{n} \leq \frac{\mathbb{E}MST_{n^2}}{n} + \frac{\mathbb{E}D_n}{n}\]
and
\[\frac{\mathbb{E}MST_k}{\sqrt{k}} \geq \frac{\mathbb{E}MST_k }{n+1} \geq \frac{\mathbb{E}MST_{n^2}}{n+1}  - \frac{\mathbb{E} D_n}{n+1}\]
and then use the fact that~\(\frac{\mathbb{E}D_n}{n} \longrightarrow 0\) as~\(n \rightarrow \infty\) (see~(\ref{dn_as})).

In the first step in the proof of~(\ref{sub_mst}), we show that
\begin{equation}\label{km_su}
\limsup_n \frac{\mathbb{E} MST_{n^2}}{n} \leq \limsup_k \frac{\mathbb{E}MST_{k^2m^2}}{km}
\end{equation}
for any fixed integer~\(m \geq 1.\)\\
\emph{Proof of~(\ref{km_su})}: Fix an integer~\(m \geq 1\) and write~\(n = qm + s\)
where~\(q = q(n) \geq 1\) and~\(0 \leq s = s(n) \leq m-1\) are integers.
As~\(n \rightarrow \infty,\)
\begin{equation}\label{kn_def}
q(n) \longrightarrow \infty \text{ and }\frac{n}{q(n)} \longrightarrow m.
\end{equation}

Using property~\((t1),\)
\[MST_{n^2} = MST_{(qm+s)^2} = MST_{q^2m^2 + 2qms + s^2} \leq MST_{q^2m^2} + 3\sqrt{2qms + s^2} + \sqrt{2}\]
and so
\begin{equation}\label{mst_temp_km1}
\limsup_n \frac{\mathbb{E}MST_{n^2}}{n} \leq \limsup_n \frac{qm}{n} \frac{\mathbb{E} MST_{q^2m^2}}{qm} + \limsup_n \frac{3\sqrt{2qms + s^2}+\sqrt{2}}{qm}.
\end{equation}
Since~\(s \leq m-1 < m,\)~\(3\sqrt{2qms + s^2} + \sqrt{2} \leq 4m\sqrt{2q+1} + \sqrt{2}\)
and so using~(\ref{kn_def}), the second term in~(\ref{mst_temp_km1}) is zero.
Using~(\ref{kn_def}) again, the first term in~(\ref{mst_temp_km1}) equals
\begin{equation}\label{su_li}
\limsup_n \frac{\mathbb{E} MST_{q^2 m^2}}{qm} \leq \limsup_k \frac{\mathbb{E} MST_{k^2m^2}}{k}.
\end{equation}
This proves~(\ref{km_su}).

\emph{Proof of~(\ref{su_li})}: Let~\(L_1 = \limsup_n \frac{\mathbb{E} MST_{q^2 m^2}}{qm}\) and~\(L_2 = \limsup_k \frac{\mathbb{E} MST_{k^2m^2}}{k}.\)
For~\(q = q(n)\) as defined prior to~(\ref{kn_def}) and for all integers~\(l \geq 1,\) we have
\[\sup_{n \geq l} \frac{\mathbb{E} MST_{q^2 m^2}}{qm} \geq L_1 \text{ and so }\sup_{ n \geq lm + m} \frac{\mathbb{E} MST_{q^2 m^2}}{qm} \geq L_1.\]
But~\(n \geq lm +m \) implies that~\(q(n) = \frac{lm+m-s}{m} \geq l\)
since~\(s = s(n) \leq m\) (see statement prior to~(\ref{kn_def})). Therefore
\[L_1 \leq \sup_{n \geq lm +m } \frac{\mathbb{E} MST_{q^2 m^2}}{qm} \leq \sup_{ k \geq l} \frac{\mathbb{E} MST_{k^2 m^2}}{km} \downarrow L_2\]
as~\(l \rightarrow \infty.\) ~\(\qed\)

%In particular if~\(l \geq 1\) is an integer then \[\sup_{ n \geq lm + m} \frac{\mathbb{E} MST_{q^2 m^2}}{qm} \geq L_1.\]

If~\(\lambda := \liminf_n \frac{\mathbb{E}MST_{n^2}}{n}\) then from~(\ref{exp_mst_u}), we have that~\(\lambda > 0.\) Moreover,
\begin{equation}\label{second_step_est}
\limsup_{k} \frac{\mathbb{E} MST_{k^2m^2}}{km} \leq \lambda
\end{equation}
and so~(\ref{beta_conv}) follows from~(\ref{km_su}).

To prove~(\ref{second_step_est}), we proceed as follows.
For positive integers~\(k\) and~\(m,\) distribute~\(k^2 m^2\) nodes~\(\{X_i\}_{1 \leq i \leq k^2m^2}\) independently
and uniformly in the unit square~\(S.\) Divide~\(S\) into~\(k^2\) disjoint
squares~\(\{W_j\}_{1 \leq j \leq k^2m^2}\) each of size~\(\frac{1}{k} \times \frac{1}{k}\)
and let~
\begin{equation}\label{nj_w_def}
N_j = \sum_{i=1}^{k^2m^2} \ind(X_i \in W_j)
\end{equation}
denote the number of nodes in the square~\(W_j.\)\\
\((t2)\) If~\(MST(N_j)\) denotes the length MST of the nodes in the square~\(W_j\)
then
\begin{equation}\label{mst_km}
MST_{k^2m^2} \leq \sum_{j=1}^{k^2} MST(N_j) + 4k\sqrt{2}.
\end{equation}
\emph{Proof of~\((t2)\)}:
For the proof of~(\ref{mst_km}), we proceed as in the proof of the strips method (see~(\ref{mst_ab})).
Suppose the top left most square is labelled~\(W_1,\) the square below~\(W_1\) is~\(W_2\) and so on until
we reach the square~\(W_k\) intersecting the bottom edge of the unit square~\(S.\) The square to the right of~\(W_1\)
is then labelled~\(W_{k+1}\) and the square below~\(W_{k+1}\) is~\(W_{k+2}\) and so on.
For~\(j \geq 1,\) let~\({\cal T}(j)\)
be the MST formed by the nodes of~\(W_j.\) We set~\({\cal T}(j) = \emptyset \)
if~\(W_j\) contains no node. Suppose~\({\cal T}(1) \neq \emptyset\)
and let~\(W_{j_1}\) be the ``first" square below~\(W_1\) in the first column of squares~\(\{W_i\}_{1 \leq i \leq k}\)
also containing at least one node. %as shown in Figure~\ref{fig_tree_join}. Here the small square
%containing the node~\(A\) is~\(W_1\) and the small square containing the node~\(B\) in~\(W_{j_1}.\)

\begin{figure}[tbp]
\centering
%\fbox{
\includegraphics[width=2in, trim= 50 180 50 100, clip=true]{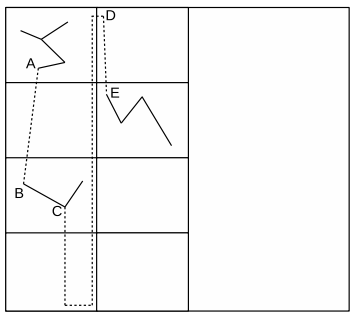}
%}
\caption{Joining trees in each subsquare to form a large spanning tree.}
\label{fig_tree_join}
\end{figure}

Join some node of~\(A \in {\cal T}(1)\) with some node~\(B \in {\cal T}(j_1)\) and call the resulting
edge as an inclined \emph{extra edge} (see Figure~\ref{fig_tree_join}). Similarly let~\(j_2 \geq j_1+1\) be the least
indexed square containing at least one node in the first column of squares~\(\{W_i\}_{1 \leq i \leq k}\)
and join some node of~\({\cal T}(j_1)\) with some node of~\({\cal T}(j_2)\)
by an inclined extra edge.

Let~\(j_{last}\) be the ``last" square in~\(\{W_i\}_{1 \leq i\leq k}\) containing at least one node
and let~\(i_1\) be the ``first" square in the second column of squares containing at least one node.
Join some node~\(C \in {\cal T}(j_{last})\) with some point~\(D\) within the first square in the second column~\(W_{k+1}\)
by vertical and horizontal extra edges as shown in Figure~\ref{fig_tree_join}. Join~\(D\) to some node of~\(E \in {\cal T}(i_1)\) by an
inclined extra edge as shown in Figure~\ref{fig_tree_join}.

Continue the above procedure for the second column of squares and proceeding iteratively, we finally obtain
a spanning tree containing all the~\(k^2m^2\) vertices. By construction, any extra edge (horizontal, vertical or inclined)
intersecting a square~\(W_j\) has length no more than~\(\frac{\sqrt{2}}{k},\) the length of the diagonal of~\(W_j.\) Also, at most
four extra edges intersect~\(W_j.\)
Since there are~\(k^2\) squares in~\(\{W_j\},\) the total length of the extra edges added is no more
than~\(4\frac{\sqrt{2}}{k}k^2 = 4k\sqrt{2}.\)~\(\qed\)

From property~\((t2)\) and~(\ref{mst_ab2}), we get
\begin{equation}\label{jtemp}
\mathbb{E} MST_{k^2m^2} \leq k^2 \mathbb{E}MST(N_1) + 4k \sqrt{2}
\end{equation}
To evaluate~\(MST(N_1),\) write~\(\mathbb{E} MST(N_1) = I_1 + I_2,\) where~\[I_1 =\mathbb{E} MST(N(1)) \ind(F_1), I_2 = \mathbb{E}MST(N(1)) \ind(F_1^c)\]
and~\(F_1 := \{m^2 - m\log{m} \leq N_1 \leq m^2 + m\log{m}\}.\)
Each node~\(X_i, 1 \leq i \leq k^2m^2\) has a probability~\(\frac{1}{k^2}\) of being present in the~\(\frac{1}{k} \times \frac{1}{k}\) square~\(W_1.\) Therefore the number of nodes~\(N_1\) in the square~\(W_1\) is binomially distributed with mean~\(\mathbb{E}N_1 = m^2\) and~\(var(N_1) \leq k^2m^2 \frac{1}{k^2} = m^2\) (see~(\ref{nj_w_def})). We therefore get from Chebychev's inequality that
\begin{equation}\label{fm_k_est}
\mathbb{P}(F^c_1)  \leq \frac{1}{(\log{m})^2} \leq \epsilon
\end{equation}
for all~\(m \geq M_0\) large, not depending on~\(k.\)

We evaluate~\(I_1\) and~\(I_2\) separately below.\\
\underline{\emph{Evaluation of~\(I_1\)}}:
Write~\(I_1 = \sum_{j = j_{low}}^{j_{up}} \mathbb{E} MST(N(1)) \ind(N(1) = j),\)
where~\(j_{low} := m^2 - m\log{m} \leq m^2 + m\log{m} =: j_{up}.\) Given~\(N_1 = j,\) the nodes in~\(W_1\) are uniformly distributed in~\(W_1\) and recall from discussion prior to~(\ref{mst_scale}) that~\(\mathbb{E}MST\left(j;\frac{1}{k}\right)\) is the
expected length of the MST containing~\(j\) nodes uniformly distributed in the~\(\frac{1}{k} \times \frac{1}{k}\) square~\(W_1.\)
Thus
\begin{equation}
I_1 = \sum_{j =j_{low}}^{j_{up}} \mathbb{E}MST\left(j;\frac{1}{k}\right) \mathbb{P}(N(1) = j) = \frac{1}{k} \sum_{j=j_{low}}^{j_{up}} \left(\mathbb{E}MST_j \right)\mathbb{P}(N(1) = j), \label{temp1}
\end{equation}
by~(\ref{mst_scale}).

Using the difference estimate~(\ref{mst_diff_one}) from Lemma~\ref{one_point_est}, we have for any
\(j_{low} \leq j_1, j_2 \leq j_{up}\) that
\[\mathbb{E}|MST_{j_2} - MST_{j_1}|  \leq \sum_{u=j_{low}}^{j_{up}-1} \mathbb{E}|MST_{u+1} - MST_{u}|
\leq \sum_{u=j_{low}}^{j_{up}-1} C\frac{(\log{u})^{3/2}}{\sqrt{u}}\]
for some constant~\(C > 0\) not depending on~\(j_1\) or~\(j_2.\)
For all~\(j_{low} \leq u \leq j_{up},\) the term~\(\frac{(\log{u})^{3/2}}{\sqrt{u}} \leq C_1 \frac{(\log{m})^{3/2}}{m}\)
for some positive constant~\(C_1\) and so the term~\(\mathbb{E}|MST_{j_2} - MST_{j_1}|  \) is bounded above by
\begin{equation}
(j_{up}-j_{low}) C_1\frac{(\log{m})^{3/2}}{m} \leq (2m\log{m}) C_1\frac{(\log{m})^{3/2}}{m}  = C_2 (\log{m})^{5/2} \label{temp3}
\end{equation}
for some constant~\(C_2 > 0.\)
Setting~\(j_1= m^2\) and~\(j_2 = j\) and using~(\ref{temp3}) we get~\(MST_j \leq MST_{m^2} + C_2 (\log{m})^{5/2}\) for all~\(j_{low} \leq j \leq j_{up}.\)
From~(\ref{temp1}) we therefore have that
\begin{equation}
I_1 \leq \frac{1}{k} \mathbb{E} MST_{m^2} + \frac{1}{k} C_2 (\log{m})^{5/2}. \label{temp6}
\end{equation}
%since there are at most~\(2\log{m}\) terms in the summation term in the first line of~(\ref{temp6}) (see~(\ref{j_range})).

\underline{\emph{Evaluation of~\(I_2\)}}: There are~\(N(1)\) nodes
in the square~\(W_1\) and so from the strips estimate~(\ref{mst_ab}),~\(MST(N(1)) \leq \frac{3}{k}\sqrt{N(1)}.\)
Thus
\begin{equation}\label{i2_estw}
I_2 = \mathbb{E}MST(N(1))\ind(F_1^c) \leq \frac{3}{k} \mathbb{E}\sqrt{N(1)} \ind(F_1^c) \leq \frac{3}{k} \left(\mathbb{E}N_1\right)^{\frac{1}{2}} \left(\mathbb{P}(F_1^c)\right)^{\frac{1}{2}},
\end{equation}
by the Cauchy-Schwarz inequality. Since~\(\mathbb{E}N_1= m^2\) and~\(\mathbb{P}(F_1^c) \leq \epsilon\)
for a fixed~\(\epsilon> 0\) and for all~\(m\) large (see~(\ref{fm_k_est})), we get
\begin{equation}\label{i2_estww}
I_2 \leq \frac{3}{k} m \sqrt{\epsilon}.
\end{equation}

Substituting~(\ref{i2_estww}) and~(\ref{temp6}) into~(\ref{jtemp}) gives
\begin{equation}\label{jtemp2}
\mathbb{E} MST_{k^2m^2} \leq k\mathbb{E} MST_{m^2} + C_2 k (\log{m})^{5/2} + 3mk\sqrt{\epsilon} + 4k\sqrt{2}
\end{equation}
and so
\begin{equation}\label{jtemp3}
\limsup_k \frac{\mathbb{E} MST_{k^2m^2}}{km} \leq \frac{1}{m}\mathbb{E} MST_{m^2} + C_2 \frac{(\log{m})^{5/2}}{m} + 3 \sqrt{\epsilon} + \frac{4\sqrt{2}}{m}
\end{equation}
for all~\(m\) large. Consequently,~\(\limsup_k \frac{\mathbb{E} MST_{k^2m^2}}{km} \leq \lambda + 3\sqrt{\epsilon}\) and since~\(\epsilon  >0\) is arbitrary, we get~(\ref{second_step_est}).\(\qed\)

%Divide the unit square into small disjoint subsquares~\(\{S_l\}\) of size~\(r_n \times r_n,\) where~\(r_n^2 = \frac{M\log{n} + a_n}{n},\) where~\(a_n \in (1,2)\) is such that~\(\frac{1}{10r_n}\) is an integer and~\(M > 0\) is a large constant to be determined later. Fix an integer~\(K \geq 1\) to be determined later and for~\(1 \leq i \leq K,\) label the centres of squares in~\(\{S_l\}\) intersecting~\((1-(i-1)100r_n)S\) as~\(\{A^{(i)}_j\}\) as shown in Figure~\ref{}. Further denote~\(S^{(i)}_j\) to be the square in~\(\{S_l\}\) with centre~\(A^{(i)}_j.\) Thus the squares~\(\{S^{(1)}_j\}\) intersect the boundary of~\(S,\) the squares~\(\{S^{(2)}_j\}\) intersect the boundary of~\((1-100r_n)S\) and so on. The squares~\(S^{(i)}_j\) are placed in such a way that the distance~\(Q_i\) between~\(A^{(i)}_j\) and~\(A^{(i)}_{j+1}\) lies between~\(10 r_n\) and~\(15 r_n.\)

%\emph{Proof of~\((a1)-(a2)\)}:
%WRTE MORE HERE +eTC...

%We now construct an approximation~\({\cal U}_n\) of~\({\cal P}_n\) and determine the number of edges in~\({\cal U}_n.\) We write...3332

\subsection*{Acknowledgement}
I thank Professors Rahul Roy, Jacob van den Berg, Anish Sarkar and Federico Camia for crucial comments and for my fellowships.

\bibliographystyle{plain}

\end{document}